\documentclass[letterpaper,11pt, nopreprintline,author year,numbers,sort&compress]{elsarticle}

\usepackage{graphicx}
\usepackage[scriptsize]{subfigure}
\usepackage{color}
\usepackage{tikz}
\usepackage{newlfont}
\usepackage{multirow}
\usepackage{array}
\usepackage[font=footnotesize]{caption} 
\DeclareCaptionStyle{period-newline}%
[labelsep=period,justification=centering]%
{labelsep=period-newline}
\usepackage{longtable} 
\usepackage{ifthen}
\usepackage{alltt}
\usepackage{float}
\usepackage{enumerate}
\usepackage[text={6.5in,8.5in},centering]{geometry} 
\newcolumntype{C}[1]{>{\centering\let\newline\\\arraybackslash\hspace{0pt}}m{#1}}
\newcolumntype{L}[1]{>{\raggedright\let\newline\\\arraybackslash\hspace{0pt}}m{#1}}
\newcolumntype{R}[1]{>{\raggedleft\let\newline\\\arraybackslash\hspace{0pt}}m{#1}}
\newcommand{\ba}{\begin{array} }
\newcommand{\ea}{\end{array} }
\newcommand{\bae}{\begin{eqnarray}}
\newcommand{\eae}{\end{eqnarray}}
\newcommand{\bea}{\begin{eqnarray*}}
\newcommand{\eea}{\end{eqnarray*}}
\newcommand{\be}{\begin{equation}}
\newcommand{\ee}{\end{equation}}

\newcommand{\lx}{\left(}
\newcommand{\rx}{\right)}
\newcommand{\lz}{\left[ }
\newcommand{\rz}{\right] }

\newcommand{\pr}{{\bf Proof}~~}

\def\qed{\mbox{}\hfill\raisebox{-2pt}{\rule{5.6pt}{8pt}\rule{4pt}{0pt}}\medskip\par}

\usepackage{amsfonts}
\usepackage{amssymb}
\usepackage{amsthm}
\usepackage{ mathrsfs }
 \usepackage{ae} 
\usepackage[T1]{fontenc}
\usepackage[ansinew]{inputenc}
\usepackage{amsmath}
\usepackage{eucal}

\usepackage[english]{babel}

\usepackage{color}
\usepackage{hyperref}
\hypersetup{
     colorlinks   = true,
     citecolor    = blue
}
\usepackage{lscape}
\usepackage{graphicx}
\usepackage{epstopdf}
\newtheorem{theorem}{\hskip\parindent\bf Theorem}[section]

\newtheorem{proposition}{\bf Proposition}[section]

\usepackage{booktabs}
\usepackage{bm}

\makeatletter
\ifcase \@ptsize \relax
  \newcommand{\miniscule}{\@setfontsize\miniscule{3}{4}}
\or
  \newcommand{\miniscule}{\@setfontsize\miniscule{4}{5}}
\or
  \newcommand{\miniscule}{\@setfontsize\miniscule{4}{5}}
\fi
\makeatother
\usepackage{lipsum}

\usepackage{float}

\usepackage{titlesec}
\titleformat*{\section}{\large\bfseries}
\titleformat*{\subsection}{\large\bfseries}

\usepackage{relsize}

\usepackage{pbox} 

\usepackage{pifont}
\newcommand{\xmark}{\ding{55}}%

\begin{document}

\begin{frontmatter}

\title{A two patch prey-predator model with multiple foraging strategies in predators: Applications to Insects}
\author[1]{Komi Messan} \ead{kmessan@asu.edu}
\address[1]{Simon A. Levin Mathematical and Computational Modeling Sciences Center, Arizona State University, \\Tempe, AZ 85281, USA.}
\address[2]{Sciences and Mathematics Faculty, College of Letters and Sciences, Arizona State University, \\Mesa, AZ 85212, USA.}
\author[2]{ Yun Kang } \ead{yun.kang@asu.edu}
 \date{Received: date / Accepted: date}


\begin{abstract} \\
 We propose and study a two patch Rosenzweig-MacArthur prey-predator model with immobile prey and predator using two dispersal strategies.  The first dispersal strategy is driven by the prey-predator interaction strength, and the second dispersal is prompted by the local population density of predators which is referred as the passive dispersal. The dispersal strategies using by predator are measured by the proportion of the predator population using the passive dispersal strategy which is a parameter ranging from 0 to 1. We focus on how the dispersal strategies and the related dispersal strengths affect population dynamics of prey and predator, hence generate different spatial dynamical patterns in heterogeneous environment.  We provide local and global dynamics of the proposed model. Based on our analytical and numerical analysis, interesting findings could be summarized as follow:  (1) If there is no prey in one patch, then the large value of dispersal strength and the large predator population using the passive dispersal in the other patch could drive predator extinct at least locally. However, the intermediate predator population using the passive dispersal could lead to multiple interior equilibria and potentially stabilize the dynamics; (2) The large dispersal strength in one patch may stabilize the boundary equilibrium and lead to the extinction of predator in two patches locally when predators use two dispersal strategies; (3) For symmetric patches (i.e., all the life history parameters are the same except the dispersal strengths), the large predator population using the passive dispersal can generate multiple interior attractors; (4) The dispersal strategies can stabilize the system, or destabilize the system through generating multiple interior equilibria that lead to multiple attractors; and (5) The large predator population using the passive dispersal could lead to no interior equilibrium but both prey and predator can coexist through fluctuating dynamics for almost all initial conditions. \\
 
 \end{abstract}


\end{frontmatter}


\section{Introduction}
The dispersal of an individual has consequences not only for individual fitness, but also for population dynamics and genetics, and species' distributions \citep{Gilpin1991, hanski1999habitat, clobert2001dispersal, bowler2005causes}. As the impact of dispersal on population dynamics has been increasingly recognized,  understanding the link between dispersal and population dynamics is vital for population management and for predicting how population responses to changes in the environment. For many animals and insects, the costs and  benefits of dispersal will vary in space and time, and among individuals. Thus, the  profit of the dispersal ability as a life-history strategy will vary as a result, and a plastic dispersal strategy is typically expected to respond to this variation \citep{ims2001condition, ronce2001perspectives, massot2002condition, bowler2005causes}. The varied dispersal driving forces include population density, kin selection relatedness, conspecific attraction, interspecific interactions, food availability, patch size and qualities, etc. There has been a large number of empirical studies supporting the effects of various parameters on dispersal   mechanisms and strengths \citep{bowler2005causes}. For example, the field work by \citet{kiester1974strategy} showed evidence of Iguanid lizards that encompass two or more dispersal strategies as foraging movements. \citet{kummel2013aphids} showed through their field work that the foraging behavior of Coccinellids are  governed not only by the conspecific attraction but also through the passive diffusion and retention on plants with high immobile aphids number. The main purpose of this article is to investigate the effects of the combinations of different strategies on population dynamics of a prey-predator interaction model when prey is immobile.  \\

Due to the practical difficulties associated with the field study of dispersal, theoretical studies play a particularly important role in predicting the effects of varied dispersal strategies in population dynamics  \citep{bowler2005causes}. The patchy prey-predator population models with different dispersal forms have  been proposed and studied  in a fair amount of literature. For example, the work of \citep{nguyen2012effects, namba1980density, janosi1997evolution, silva2001stability, hansson1991dispersal, fraser1982experimental, savino1989behavioural} explored the effects of dispersal on population dynamics of prey-predator models when local population density is a selecting factor for dispersal. The work of \citet{huang2001predator} and \citet{ghosh2011two} studied the population dynamics of a two patch model with dispersal in predator driving by local population density of prey through Holling searching-handling time budget argument. The work of \citet{kareiva1987swarms} studied dynamics when the dispersal of predator is carried out due to the concentrated food resources.  \citet{cressman2013two} investigated a two patch population-dispersal dynamics for predator-prey interactions with dispersal directed by the fitness. Recent work of \citep{kang2014dispersal} studied a two patch prey-predator model where predator is dispersed to the patch with the stronger strength of prey-predator interaction. These theoretical work provide useful insights on the link of dispersal strategies and prey-predator population dynamics.\\

Many empirical work of animal and insects show that  dispersal strategies vary among species according to their life history and how they interact with the environment \citep{bowler2005causes}. However,  there is a limited theoretical work on studying how combinations of different dispersal strategies affect population dynamics of prey-predator models in the patchy environment. This paper presents an extended version of a Rosenzweig-MacArthur two patch prey-predator model studied in \citep{kang2014dispersal} where prey is immobile and the dispersal of predator is attracted by the strength of prey-predation interaction.  Our proposed model is motivated by the field experiments of \citep{kummel2013aphids, stamps1988conspecific, kiester1974strategy}. The current model integrates the two dispersal strategies of predator: (1) the passive dispersal, i.e., the classical foraging behavior where predator is driven to the patch with the lower predator population density (e.g. \citep{jansen1995regulation}); (2) the density dependent dispersal measured through the predation attraction  \citep{kang2014dispersal}. The linear combination of these two strategies is linked through a parameter whose value is between 0 and 1, and measures the proportion of the predator population using these two dispersal strategies.  We aim to use our model to explore how the combinations of these two dispersal strategies of predator affect population dynamics of prey-predator interaction. \\

The paper is organized as follows: Section 2 introduces the proposed model along with its biological derivation, and provide a brief summary on the dynamics of the related subsystems. Section 3 presents mathematical analysis of the local and global dynamics of the proposed model. Section 4 Investigates the effects of dispersal strategies through bifurcation diagrams. Section 5 concludes our findings along with the related potential biological interpretations.\\

\section{Model derivations and the related dynamics}\label{ModelDerivation}
Let $x_i(t),~y_i(t)$ be the population of prey and predator in Patch $ i$ at time $t$, respectively. In the absence of dispersal, we assume that the population dynamics of prey and predator follow the Rosenzweig-MacArthur prey-predator model.  The dispersal of predator from Patch $i$ to Patch $j$ is driven by two mechanisms. The first mechanism relies on the strength of the prey-predation interaction in Patch $j$ (also called ``the predation strength"). Let $\rho_i$ represents the relative dispersal rate of predator at Patch $i$, then we obtained the following net predation attraction driven dispersal of predator at Patch $i$

$$\rho_i\left(\frac{a_j x_j y_j}{1 + x_j} y_i-\frac{a_i x_i y_i}{1 + x_i} y_j\right).$$ 
This assumption follows directly from the experimental work of \citet{stamps1988conspecific} in which he concluded that Anolis aeneus juveniles are attracted to conspecific territorial residents under natural conditions in the field.  This assumption has also been supported by many field studies including \citep{hassell1978foraging, alonso2004distribution, auger1996emergence}. \\

The second dispersal mechanism is termed as ``the passive dispersal" in which the dispersal is driven by the local population density of predator. The effects of this dispersal strategy has been well studied by many researchers \citep{jansen1995regulation, matthysen2005density, nguyen2012effects, poggiale1998behavioural, namba1980density, janosi1997evolution, silva2001stability, hastings1983can}. For example overcrowding of predator in a patch may decrease the resource assessment that can constitute a cue for for the local predators to move. Following this inference, the net dispersal of predators from Patch $i$ to Patch $j$ is given by 

$$\rho_i(y_j-y_i).\\$$

Motivated by the field work of \citet{kiester1974strategy} on Iguanid lizards and \citet{kummel2013aphids} on Coccinellids, we incorporate these two dispersal strategies above into our model. After similar rescaling approach by \citet{liu2003complex},  our proposed model is presented as follows with $r_1=1, r_2 = r$ being the relative intrinsic growth rates, $K_i$ being the relative carrying capacity of prey at Patch $i$ in the absence of predation, $d_i$ being the death rate of predator in Patch $i$, and the parameter $ s \in [0,1]$ representing the proportion of predator population using the passive dispersal strategy:\\

{\footnotesize
\begin{equation}
\begin{aligned}
\frac{dx_1}{dt}  &= r_1x_1\left(1 - \frac{x_1}{K_1}\right) - \frac{a_1 x_1 y_1}{1 + x_1} \\
\frac{dy_1}{dt}  &=  \frac{a_1 x_1 y_1}{1 + x_1} -d_1 y_1+ \rho_1(1-s)\left(  \underbrace{ \frac{a_1 x_1 y_1}{1 + x_1}}_\text{attraction strength to Patch 1} y_2- \underbrace{\frac{a_2 x_2 y_2}{1 + x_2}}_\text{attraction strength to Patch 2} y_1\right)+ \rho_1s(y_2-y_1)\\
\frac{dx_2}{dt} & = r_2 x_2\left(1 - \frac{x_2}{K_2}\right) - \frac{a_2 x_2 y_2}{1 + x_2} \\
\frac{dy_2}{dt} & =  \frac{a_2 x_2 y_2}{1 + x_2} -d_2 y_2+\rho_2(1-s)\left(   \underbrace{\frac{a_2 x_2 y_2}{1 + x_2}}_\text{attraction strength to Patch 2} y_1-  \underbrace{\frac{a_1 x_1 y_1}{1 + x_1}}_\text{attraction strength to Patch 1} y_2\right) + \rho_2s(y_1-y_2) 
 \label{2DPatch}
\end{aligned}
\end{equation}
} First, we have the following theorem regarding the basic dynamic properties of Model \eqref{2DPatch}:\\

\begin{theorem}\label{th1:pb} 
Assume  that all parameters are positive. Model \eqref{2DPatch} is positively invariant and bounded in $\mathbb R^4_+$.  In addition, the set $\{(x_1,y_1,x_2,y_2)\in\mathbb R^4_+:x_i=0 \}$ is invariant for both $i=1,2$.\\
\end{theorem}

Our main focus is to explore how the combination of two different dispersal strategies measured by the parameter $s\in [0,1]$ affect the two patch population dynamics. Before we continue, we first provide a summary of the dynamics of the subsystems of Model \eqref{2DPatch} including the cases of $s=0$ and $s=1$. \\

 For convenience,  let $\mu_i=\frac{d_i}{a_i-d_i}, \mbox{  and   } \nu_i=\frac{r_i(K_i-\mu_i)(1+\mu_i)}{a_iK_i}$ $i=1,2$, then in the absence of  dispersal in predator, Model \eqref{2DPatch} is reduced to the following \citet{rosenzweig1963graphical} prey-predator single patch models $i=1,2$
\begin{equation}
\begin{aligned}
\frac{dx_i}{dt} & =r_i x_i\left(1 - \frac{x_i}{K_i}\right) - \frac{a_i x_i y_i}{1 + x_i} \\
\frac{dy_i}{dt} & =  \frac{a_i x_i y_i}{1 + x_i} -d_i y_i \label{Onepatch}
\end{aligned}
\end{equation} with $r_1=1$ and $r_2=r$ and its global dynamics which can be summarized from the work of \citep{liu2003complex, hsu1977mathematical, hsu1978global} as follows:
\begin{enumerate}
\item Model \eqref{Onepatch} always has two boundary equilibria $(0, 0)$, $(K_i,0)$ where the extinction $(0, 0)$ is always a saddle.
\item The boundary equilibria $(K_i,0)$ is globally asymptotically stable if $\mu_i > K_i$.
\item If $\frac{K_i-1}{2} < \mu_i < K_i$, then $(K_i,0)$ becomes saddle and the unique interior equilibria $(\mu_i,\nu_i)$ emerges which is globally asymptotically stable.
\item If $0 < \mu_i < \frac{K_i-1}{2}$, the boundary equilibrium $(K_i,0)$ is a saddle, and the unique interior equilibrium $(\mu_i,\nu_i)$ is a source where Hopf bifurcation occurs at $\mu_i = \frac{K_i - 1}{2}$. The system \eqref{Onepatch} has a unique stable limit cycle.\\
\end{enumerate}

The summary on the dynamics of Model \eqref{2DPatch} when the dispersal of  predator foraging activities is  driven by local population density (i.e., $s=1$) and when the dispersal of  predator foraging activities is  driven by predation strength (i.e. $s=0$) are briefly presented in Table \ref{table_stability} (see \cite{kang2014dispersal} for more detailed summary on the global dynamics). \\


\section{Mathematical analysis}\label{SecMathAnalysis}

From Theorem \ref{th1:pb}, we know that the set $\{(x_1,y_1,x_2,y_2)\in\mathbb R^4_+:x_i=0 \},$ is invariant for both $i=1,2$.  Assume that $x_j=0$, Model \eqref{2DPatch} is reduced to the following three species subsystem:

\begin{equation}
\begin{aligned}
\frac{dx_i}{dt} & =r_i x_i\left(1 - \frac{x_i}{K_i}\right) - \frac{a_i x_i y_i}{1 + x_i} \\
\frac{dy_i}{dt} & =  \frac{a_i x_i y_i}{1 + x_i} -d_i y_i+\rho_i(1-s)\left(\frac{a_ix_iy_i}{1+x_i}y_j\right)+ \rho_is( y_j-y_i)\\
\frac{dy_j}{dt} & = -d_j y_j-\rho_j(1-s)\left(\frac{a_ix_iy_i}{1+x_i}y_j\right)- \rho_js( y_j-y_i)
\label{xi0}
\end{aligned}
\end{equation} 
whose basic dynamics  are provided in the following theorem:\\

\begin{theorem}\label{th4:persist}[Basic dynamics of Model \eqref{xi0}]
Let  $\mu_i=\frac{d_i}{a_i-d_i}$ and $s\in (0,1)$, then the following statements of Model \eqref{xi0} are held:
\begin{enumerate}  
\item Prey $x_i$ is persistent with $\limsup_{t\rightarrow\infty} x_i(t)\leq K_i$.
\item If $\mu_i > K_i$, then predators in two patches go extinct, and the system \eqref{xi0} has global stability at $(K_i,0,0)$.  
\item  If $\rho_is<\frac{(a_i-d_i)(K_i-\mu_i)}{1+K_i}$, then predators in the two patches are persistent.\\
 \end{enumerate}
 \end{theorem}

\noindent\textbf{Notes:} Model \eqref{xi0} can apply to the case where Patch $i$ is the source patch with prey population and Patch $j$ is the sink patch without prey population. The predator in the sink patch is migrated from the source patch. Theorem \ref{th4:persist} indicates the follows regarding the effects of the proportion of predator using the passive dispersal on Model \eqref{xi0}:
\begin{enumerate}
\item Prey $x_i$ of Model \eqref{xi0} is always persistent for all $r_i>0$. This is different than the case of $s=1$ since prey may go extinct when $s=1$.
\item If $\mu_i<K_i$ and $\rho_is$ is small enough, then  the inequality $\rho_is<\frac{(a_i-d_i)(K_i-\mu_i)}{1+K_i}$ holds, hence predators persist. This result suggests that, under the condition of $\mu_i<K_i$, the large value of $\rho_i s$ could drive predator extinction in two patches at least locally.\\
\end{enumerate}


The interior equilibria $(x_1^*,y_1^*,y_2^*)$ of Model \eqref{xi0} is determined by first solving for $y_i^*$ and $x_i^*$ in $\frac{dx_i}{dt}=0$ and $\frac{dy_j}{dt}=0$ as follows:

{\small
\begin{equation}\label{IntB1}
\begin{aligned}
\frac{dx_i}{dt} &= 0 \quad  \Rightarrow \quad y_{i}^* =\frac{r_i(K_i-x_i^*)(1+x_i^*)}{a_iK_i} \\
\frac{dy_j}{dt} &= 0 \quad \Rightarrow \quad  x_{i}^* = \frac {-\rho_jsy_i^*+\rho_jsy_j^*+d_jy_j^*}{ a_i\rho_jsy_i^*y_j^*-a_i\rho_jy_i^*y_j^*+\rho_jsy_i^*-\rho_jsy_j^*-d_jy_j^*}
\end{aligned}
\end{equation}
}

An equation of $y_j^*$ is obtained by solving the following  equation from Model \eqref{xi0}:
{\small
\begin{align}\label{IntB2}
\frac{dy_i}{dt}\rho_j+\frac{dy_j}{dt}\rho_i = \frac{a_ix_iy_i}{1+x_i}\rho_j-d_iy_i\rho_j-d_jy_j\rho_i=0 \quad \Rightarrow \quad y_j^* = \frac{\rho_jy_i^* (-a_ix_i^*+d_ix_i^*+d_i)}{d_j\rho_i (x_i^*+i)}
\end{align}
}
A substitution of $y_i^*$ from \eqref{IntB1} into $y_j^*$ gives $y_j^* = \frac{r_i(K_i-x_i^*)[x_i^*(a_i-d_i)-d_i]\rho_j}{a_iK_id_j\rho_i}$. The discussion above implies that the existence of the interior equilibrium requires $a_i>d_i$ and $\mu_i=\frac{d_i}{a_i-d_i}<x_i^*<K_i$ otherwise $y_j^*<0$ or $y_i^*<0$. Define
{\small
\begin{align*}
f_{t_i}(x_i)&=K_i(1+x_i)[(a_i-d_i)(d_j+s\rho_j)-d_js\rho_i]-K_ia_i(d_j+s\rho_j)\\
f_{b_i}(x_i)&=[d_i-(a_i-d_i)x_i][K_i(d_j+s\rho_j)+r_i\rho_j(1-s)(1+x_i)(1+K_i)]+d_jK_is(1+x_i)\rho_i.
\end{align*}
}
Then we can conclude that  $x_i^*$ solving from Equation \eqref{IntB1} is in term of $y_i^*$ and $y_j^*$.  Upon substitution of $y_i^*$ and $y_j^*$ into $x_i^*$ we obtain the following nullclines: 

{\small
\begin{equation}\label{NullBound}
\begin{aligned}
 x_i &= \frac{K_i(1+x_i)[(a_i-d_i)(d_j+s\rho_j)-d_js\rho_i]-K_ia_i(d_j+s\rho_j)}{[d_i-(a_i-d_i)x_i][K_i(d_j+s\rho_j)+r_i\rho_j(1-s)(1+x_i)(1+K_i)]+d_jK_is(1+x_i)\rho_i} =\frac{f_{t_i}(x_i)}{f_{b_i}(x_i)} \\
& \Leftrightarrow \\
&x_if_b(x_i) - f_t(x_i) = \underbrace{[x_i^{3}-(\mu_i+K_i)x_i^{2}-\alpha_i x_i+\beta_i]}_{f_i(x_i)}[x_i+1] = 0
\end{aligned}
\end{equation}
 }
 with $\beta_i = \frac{\left[d_j\rho_is+d_i(d_j+\rho_js)\right]K_i}{r_i(a_i-d_i)(1-s)\rho_j}$ and
  {\footnotesize $$\alpha_i =\frac{\left[d_js\rho_i+r_id_i(1-s)-(a_i-d_i)(d_j+s\rho_j)\right]K_i}{r_i(a_i-d_i)(1-s)\rho_j} =\beta_i +\frac{[r_id_i(1-s)-a_i(d_j+s\rho_j)]K_i}{r_i(a_i-d_i)(1-s)\rho_j} .$$}
 
 Based on  the  arguments above and additional analysis, we have the following proposition regarding the existence of the interior equilibria of Model \eqref{xi0}:\\

 
 \begin{proposition}\label{pr1:be}  Model \eqref{xi0} can have up to two interior equilibria $E_{x_i,y_i,y_j}^{\ell}=(x_{i\ell}^*,y_{i\ell}^*,y_{j\ell}^*), \ell=1,2$. More specifically, 
\begin{enumerate}
\item If  $a_i<d_i$ or $K_i<\mu_i$ or $(\mu_i+K_i)^2+3\alpha_i< 0$, Model \eqref{xi0} has no interior equilibrium.
\item If $\frac{3\beta_i}{\mu_i+K_i}<\alpha_i<(\mu_i+K_i)^2$, then $f_i(x_i)$  has two positive roots $x_{i\ell}^*, \ell=1,2.$ If, in addition, $\mu_i<x_{i\ell}^*<K_i$ for both $\ell=1,2$, then Model \eqref{xi0} has two interior equilibria.\\
\end{enumerate}
\end{proposition}


 
 \noindent\textbf{Notes:} Proposition \eqref{pr1:be} implies that even if $f_i(x_i) $ has two positive real roots, Model \eqref{xi0} may have none or one interior equilibrium unless these two positive roots are in $(\mu_i, K_i)$. Note that the interior equilibria of the subsystem Model \eqref{xi0} represent the boundary equilibria of Model \eqref{2DPatch} when $x_1 = 0 (i=2)$ or $x_2 = 0 (i=1)$. The existence of these boundary equilibria of Model \eqref{2DPatch} when $x_1 = 0$ or $x_2 = 0$ are hence guarantee by the conditions to obtain  the interior equilibria $E_{x_i,y_i,y_j}^{\ell}$ and $E_{y_j,x_i,y_i}^{\ell}$ from Proposition \eqref{pr1:be}.  \\

In order to capture the dynamics of the interior equilibria of Model \ref{xi0}, we perform bifurcation simulations with respect to the proportion of predators using the passive dispersal, i.e., the values of $s$. Our analysis implies that Model \eqref{xi0} can have up to two interior equilibria $E_{x_1,y_1,y_2}^{\ell} = (x_{1\ell}^*,y_{1\ell}^*,y_{2\ell}^*)$ (for $i=1$) and $E_{y_1,x_2,y_2}^{\ell}= (\hat{y}_{1\ell}^*,\hat{x}_{2\ell}^*,\hat{y}_{2\ell}^*)$ $\ell=1,2$ (for $i=2$). We fix the following parameter values, 
$$r_1=1,\,\,r_2=1.8,\,\,d_2=0.35,\,\,K_1=10,\,\,K_2=7,\,\,a_2=1.4,\,\,\rho_1 = 1, \,\,\rho_2 = 2.5.$$
These fixed values implies that at Patch 2, prey and predator  coexist  in the form of a unique stable limit cycle in the absence of dispersal since $\mu_2=\frac{d_2}{a_2-d_2}=35/105<(K_2-1)/2=3$.  We consider the following two typical cases regarding the population dynamics of prey and predator in the absence of dispersal:
\begin{enumerate}
 \item $d_1 = 0.85,~a_1=1$:  Predator and prey are persistent and  have global equilibrium dynamics at Patch 1 in the absence of dispersal since $(K_1-1)/2=4.5<\mu_1=\frac{d_1}{a_1-d_1}=17/3<10=K_1$.
 \item $d_1 = 2,~a_1=2.1$: Predator goes extinct globally at Patch 1 in the absence of dispersal since $\mu_1=\frac{d_1}{a_1-d_1}=20>K_1=10$.\\
\end{enumerate}

\begin{figure}[H]
\begin{center}
   \subfigure[Effect of dispersal strategy when \newline \hspace*{.45cm} $x_2=0$ and $d_1 = 0.85,~a_1 = 1$]{\includegraphics[height = 51mm, width = 54mm]{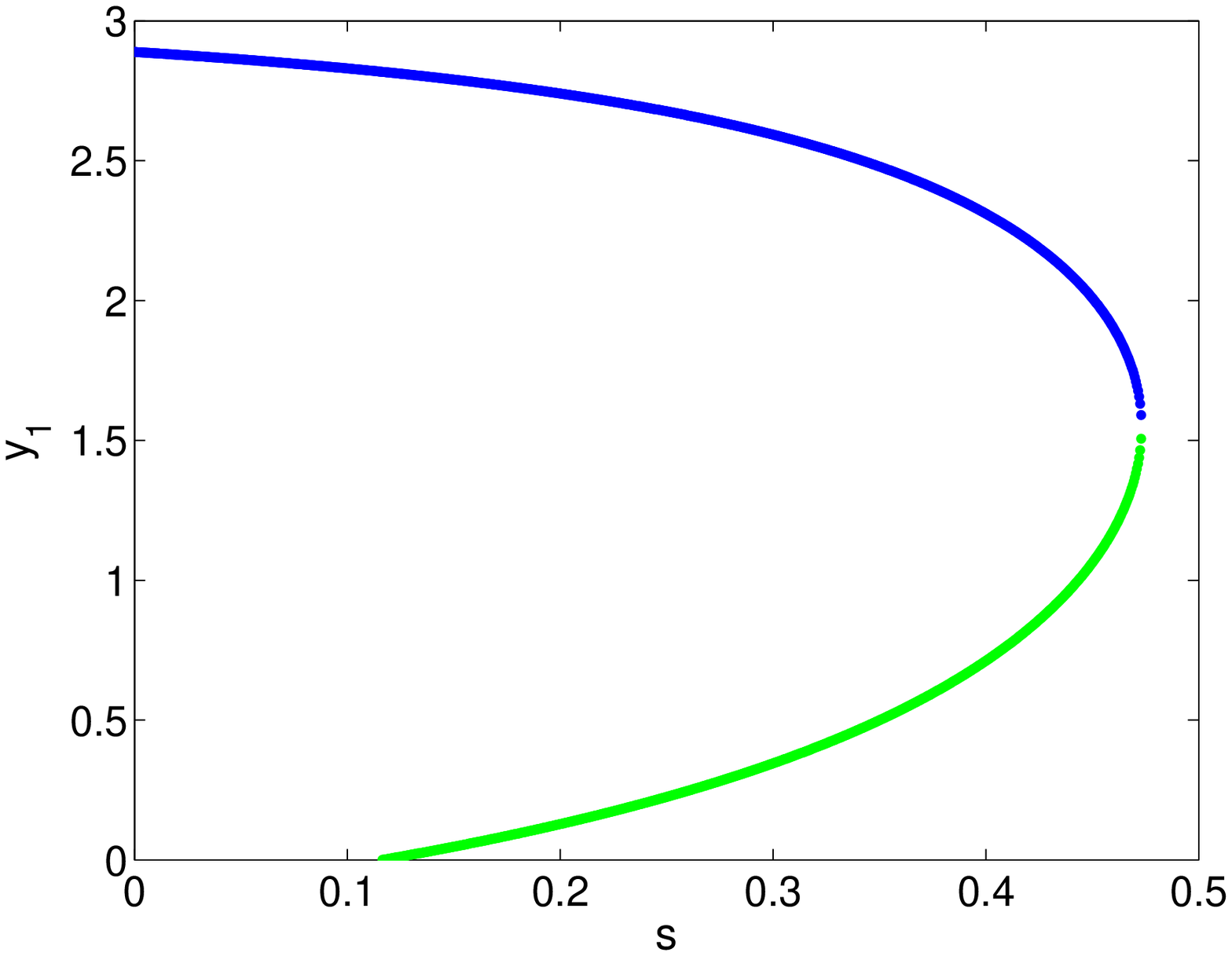}\label{fig_InteriorModel3x20_1}}
\subfigure[Effect of dispersal strategy when  \newline \hspace*{.45cm} $x_1 = 0$ and $d_1 = 0.85,~a_1 = 1$.]{\includegraphics[height = 51mm, width =54mm]{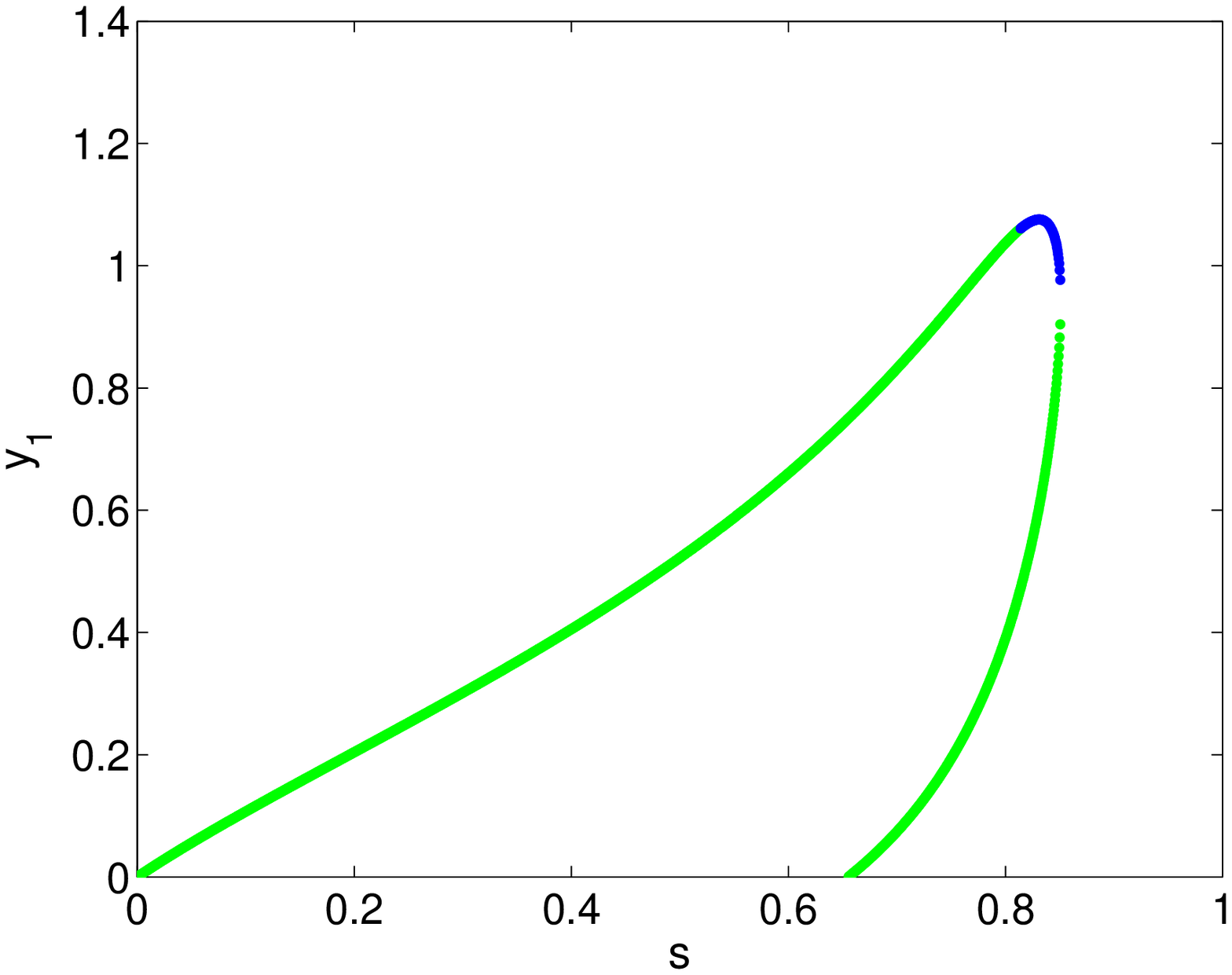}\label{fig_InteriorModel3x10_1}}\hspace{.4mm}
   \subfigure[Effect of dispersal strategy when  \newline \hspace*{.45cm} $x_1 = 0$ and $d_1 = 2,~a_1 = 2.1$.]{\includegraphics[height = 51mm, width =54mm]{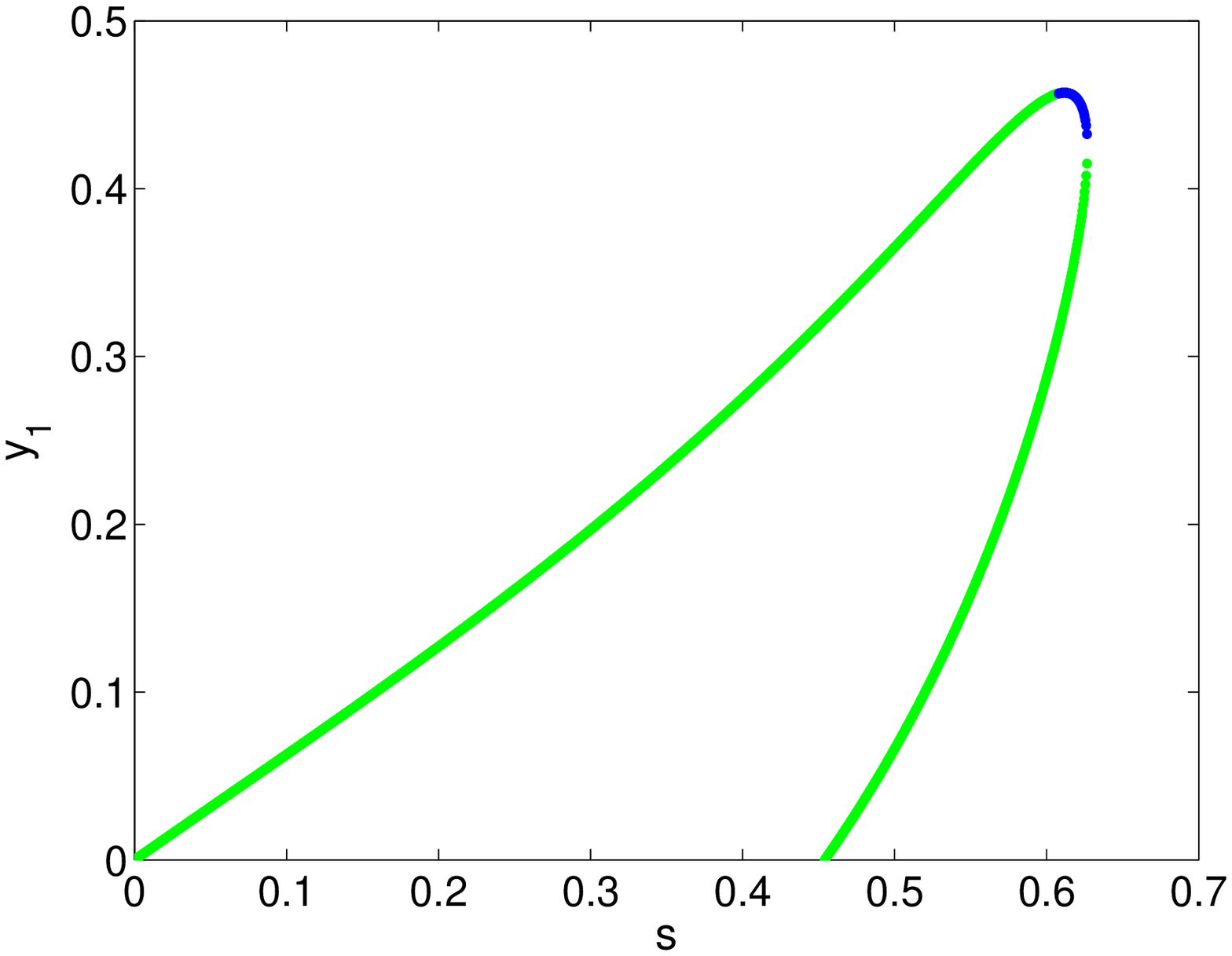}\label{fig_InteriorModel3x10_2}}\hspace{.4mm}
\end{center}
\vspace{-15pt}
\caption{One parameter bifurcation diagrams of Model \eqref{xi0} with $y$-axis representing the population size of predator at Patch 1 and $x$-axis represent the proportion of predator using the passive dispersal. Figure \ref{fig_InteriorModel3x20_1} describes  the number of interior equilibria $(\hat{x}_1^*,\hat{y}_1^*,\hat{y}_2^*)$ when $x_2 = 0$ in Model \eqref{2DPatch} and their stability with respect to variation in $s$. Figure \ref{fig_InteriorModel3x10_1} and \ref{fig_InteriorModel3x10_2} describe the number of interior equilibria $(y_1^*,x_2^*,y_2^*)$ of the submodel $x_2 = 0$ of Model \eqref{2DPatch}  and their stability  when $s$ varies from $0$ to $1$. Blue represents the sink and green represents the saddle.}
\label{fig:intm3_1}
\end{figure}

The fixed values of parameters and the two cases above provide the following four scenarios:
\begin{enumerate}
\item $i=1$ (i.e., $x_2=0$ for Model \eqref{2DPatch}) with $d_1 = 0.85,~a_1=1$. In this case, Patch 1 is the source patch and Model \eqref{xi0} can have up to two interior equilibria depending on the values of $s$ (see Figure \ref{fig_InteriorModel3x20_1}).

\item $i=1$ (i.e., $x_2=0$ for Model \eqref{2DPatch}) with $d_1 = 2,~a_1=2.1$. In this case, Patch 1 is the source patch and Model \eqref{xi0} has no interior equilibria according to Proposition \eqref{pr1:be}.

\item $i=2$ (i.e., $x_1=0$ for Model \eqref{2DPatch}) with $d_1 = 0.85,~a_1=1$. In this case, Patch 2 is the source patch and Model \eqref{xi0} can have up to two interior equilibria depending on the values of $s$ (see Figure \ref{fig_InteriorModel3x10_1}). The relative large value of $s$ can stablize the dynamics (see the blue region of Figure \ref{fig_InteriorModel3x10_1}). 

\item $i=2$ (i.e., $x_1=0$ for Model \eqref{2DPatch}) with $d_1 = 2,~a_1=2.1$. In this case, Patch 2 is the source patch and Model \eqref{xi0} can have up to two interior equilibria depending on the values of $s$ (see Figure \ref{fig_InteriorModel3x10_2}). The relative large value of $s$ can stablize the dynamics (see the blue region of Figure \ref{fig_InteriorModel3x10_2}). 
\\
\end{enumerate}

The bifurcation diagrams (Figure \ref{fig:intm3_1}) suggest that the proportion of predators using the passive dispersal can have huge impacts on the number of interior equilibria of Model \eqref{xi0}: For the small values of $s$, Model \eqref{xi0} can have one interior equilibrium ($E_{x_1,y_1,y_2}^{1}$ or $E_{y_1,x_2,y_2}^{1}$); For the intermediate values of $s$, Model \eqref{xi0} can have two interiors $E_{x_1,y_1,y_2}^{l},l=1,2$ ($i=1$) or $E_{y_1,x_2,y_2}^{l},l=1,2$ ($i=2$); For the large values of $s$, it has no interior equilibria. A more detail description of the effects of $s$ on the interior equilibria of Model \eqref{xi0} is provided in Table \eqref{table_InteriorModel_xi0}.


{\small\begin{table}[H]
\centering
\begin{tabular}{|c|C{1.4cm}|C{1.4cm}|C{1.4cm}|C{1.4cm}|C{1.4cm}|C{1.4cm}|C{1.4cm}|}
	\hline
\multicolumn{1}{|c|}{} &  \multicolumn{4}{c|}{$\mathbf{a_1= 1 \mbox{  and  } d_1 =0.85 }$ }  & \multicolumn{3}{c|}{$\mathbf{a_1= 2.1 \mbox{  and  } d_1 = 2 }$ }\\
\cline{2-8}	 
{\bf Scenarios }   &  $\mathbf{E_{x_1y_1y_2}^{1} }$    &   $\mathbf{E_{x_1y_1y_2}^{2} }$ &$\mathbf{ E_{y_1x_2y_2}^{1}}$ & $\mathbf{E_{y_1x_2y_2}^{2} }$ &  $\mathbf{E_{x_1y_1y_2}^{1,2} }$ & $\mathbf{ E_{y_1x_2y_2}^{1}}$& $\mathbf{E_{y_1x_2y_2}^{2} }$   \\ \hline
$s\leq0.1$ &  LAS &  \xmark & Saddle &  \xmark & \xmark & Saddle & \xmark \\ \hline
$ 0.15 \leq s \leq 0.45$ &   LAS & Saddle & Saddle &  \xmark & \xmark & Saddle & \xmark  \\ \hline
$ 0.55 \leq s \leq 0.62$ &  \xmark & \xmark & Saddle &  \xmark & \xmark & LAS & Saddle \\ \hline
$ 0.68 < s <  0.82$ &  \xmark &  \xmark & LAS &  Saddle & \xmark & \xmark & \xmark \\
\hline
$ s \geq 0.82$ & \xmark &  \xmark &  \xmark &  \xmark & \xmark & \xmark & \xmark \\
\hline
\end{tabular}
\caption{{\footnotesize Summary of the effect of the proportion of predators using the passive dispersal on Model \eqref{xi0} From Figures \ref{fig_InteriorModel3x20_1}, \ref{fig_InteriorModel3x10_1}, and \ref{fig_InteriorModel3x10_2}. LAS refers to local asymptotical stability and \xmark \hspace{.01in} implies the equilibrium does not exist.}}
\label{table_InteriorModel_xi0}
\end{table}}


\subsection{Boundary equilibria  and global dynamics of Model \eqref{2DPatch}}
First,  we have boundary equilibria and global dynamics of Model \eqref{2DPatch} in the following theorem.

\begin{theorem}\label{th2:be}[Boundary equilibria and global dynamics of Model \eqref{2DPatch}] Assume that $s \in (0,1)$.
Model \eqref{2DPatch} always has the following four boundary equilibrium
$$E_{0000}, E_{K_1000}, E_{00K_20}, E_{K_10K_20}$$
with the first three always being saddle.  $E_{K_10K_20}$ is locally asymptotically stable if the follwoing two inequalities in \eqref{Ek10k20st} hold:

{\footnotesize
\begin{equation}\label{Ek10k20st}
\begin{aligned}
& \mathlarger{\sum}_{i=1}^{2}\left[ \frac{(a_i-d_i)(\mu_i-K_i)}{1+K_i} + s\rho_i \right] > 0 \\
&\mbox{  and  } \\
&\left[\frac{(a_1-d_1)(\mu_1-K_1)}{1+K_1}\right]\left[s\rho_2+\frac{(a_2-d_2)(\mu_2-K_2)}{1+K_2}\right]+s\rho_1\left[\frac{(a_2-d_2)(\mu_2-K_2)}{1+K_2}\right] > 0.
\end{aligned} 
\end{equation}
}
And $E_{K_10K_20}$ is saddle when one or both of equations \eqref{Ek10k20st} are not satisfied. In addition, 
\begin{enumerate}
\item Model \eqref{2DPatch} is globally stable at $E_{K_10K_20}$ if $\mu_i>K_i$ for both $i=1,2$.
\item At least prey population in one patch of Model \eqref{2DPatch} is persistent, and the predator population in each patch is persistent if $\mu_i<K_i$ for both $i=1,2$.
\end{enumerate}
\end{theorem}
\noindent\textbf{Notes:} Theorem \ref{th2:be} indicates that the global stability of the boundary equilibrium $E_{K_10K_20}$ does not depend on the proportion of predator population using the passive dispersal since  $E_{K_10K_20}$ is globally asymptotically stable when $\mu_i>K_i, i=1,2$ which is independent of $s$. However, the value of $s>0$ and $\rho_i, i=1,2$ can stabilize $E_{K_10K_20}$. For example, assume that $\mu_i<K_i$ and $\mu_j>K_j$, then in the absence of dispersal, the boundary equilibrium $E_{K_10K_20}$ is a saddle. In the presence of the dispersal, according to Theorem \ref{th2:be},  if we choose $\rho_j$ large enough, then $E_{K_10K_20}$ can be locally stable, thus the large dispersal at one patch may stabilize the boundary equilibrium $E_{K_10K_20}$. However, if $s=0$, then dispersal has no such effects. \\

\begin{figure}[H]
\begin{center}
\subfigure[{\scriptsize Number of boundary equilibria when $x_2 = 0$}]{\includegraphics[scale=.40]{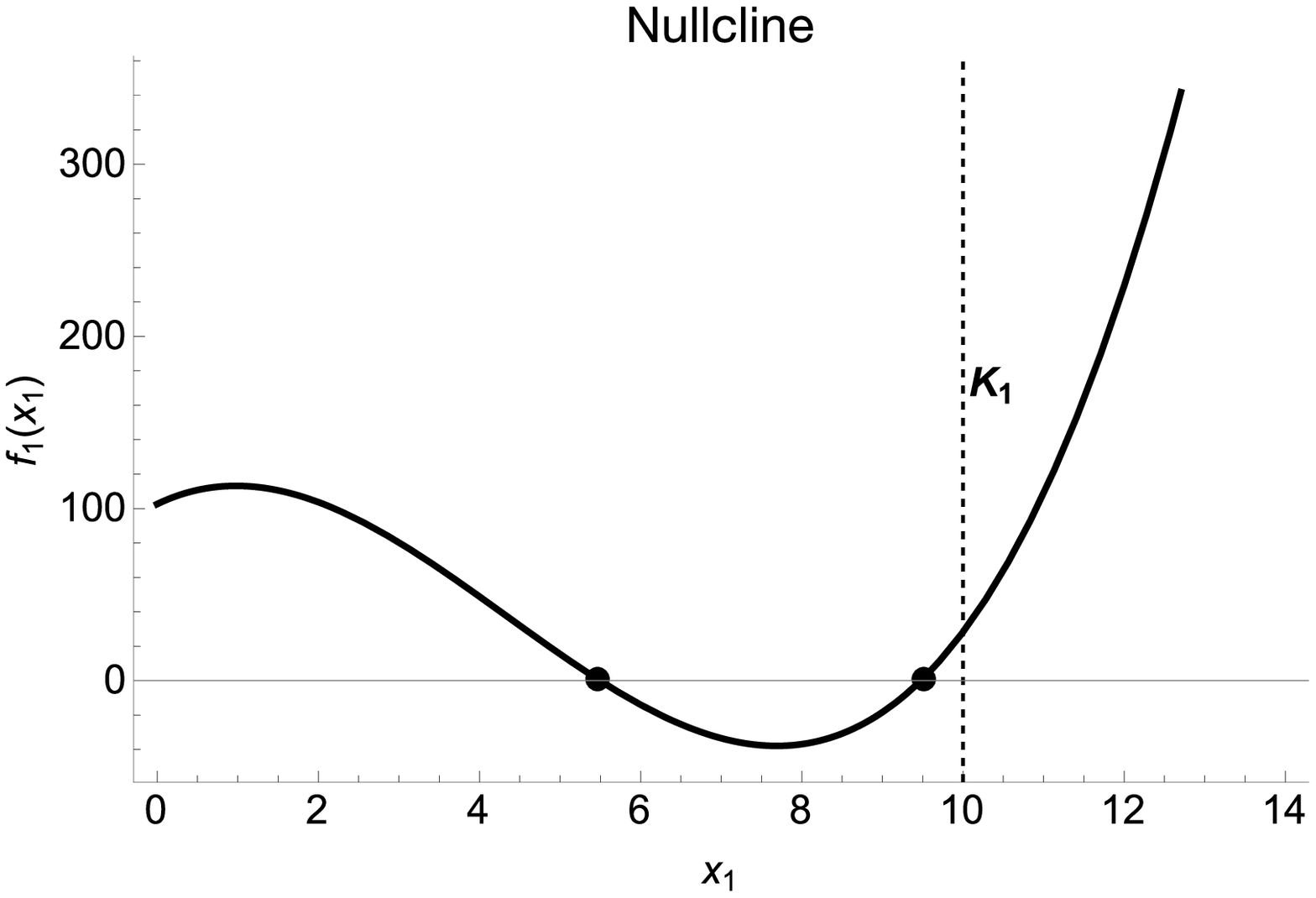}\label{fig_number_boundary_x202}}\hspace{2mm}
\subfigure[{\scriptsize Number of boundary equilibria when $x_1 = 0$}]{\includegraphics[scale=.40]{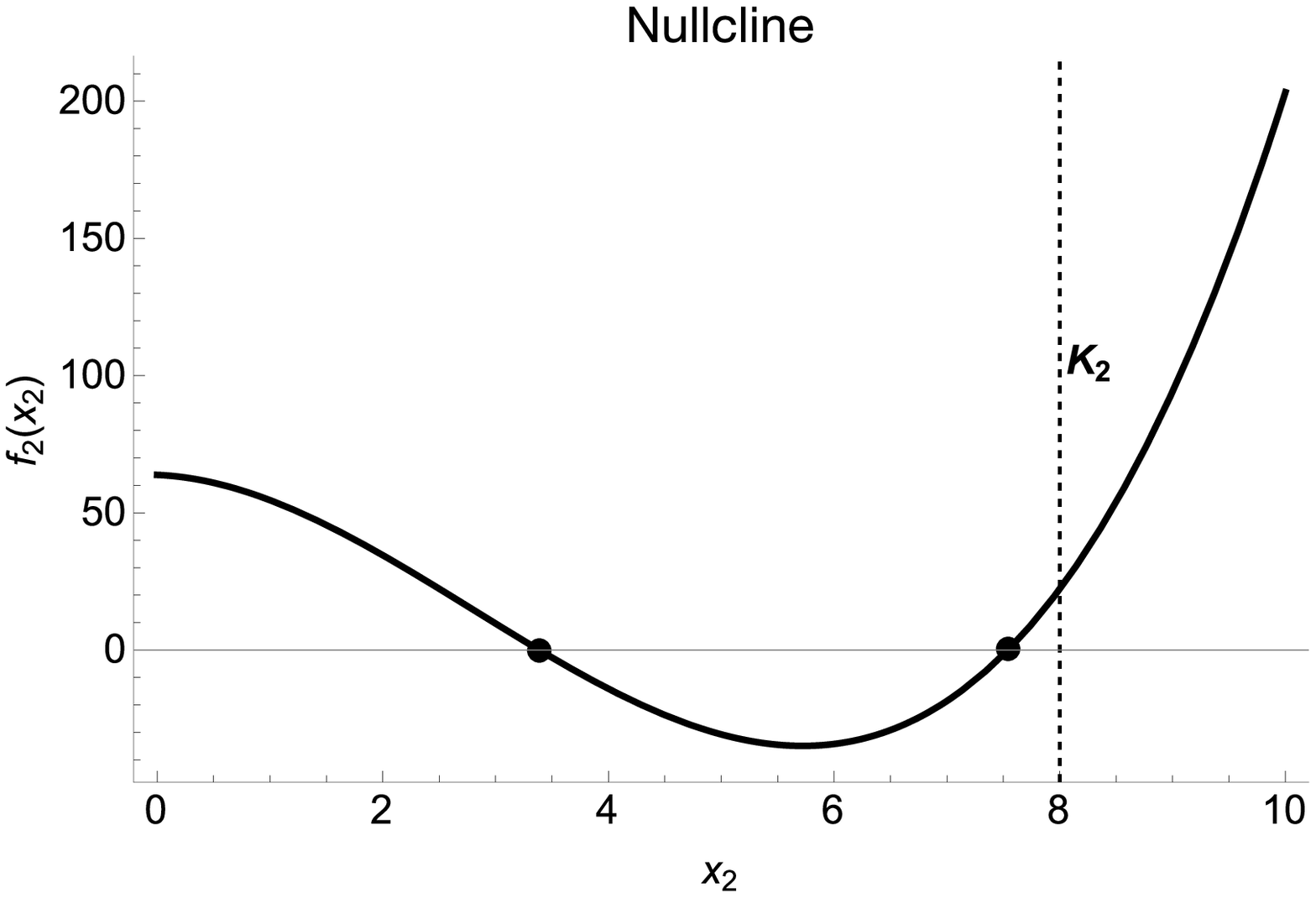}\label{fig_number_boundary_x102}}
\end{center}
\vspace{-15pt}
\caption{Boundary equilibria $E_{1\ell}^b = (x_{1\ell}^*,y_{1\ell}^*,0,y_{2\ell}^*)$ and $E_{2\ell}^b = (0,\hat{y}_{1\ell}^*,\hat{x}_{2\ell}^*,\hat{y}_{2\ell}^*)$. Figure  \ref{fig_number_boundary_x202} and \ref{fig_number_boundary_x102} are the cases when $s = 0.65$, $r_1=1$, $r_2=0.54$, $d_1=0.45$, $d_2=0.105$, $K_1=10$, $K_2=8$, $a_1=0.6$, $a_2=0.35$, $\rho_1 = 1.75$, and $\rho_2=1.2$. The solid lines are $f_1(x_1)$ and $f_2(x_2)$ while the dashed lines are $K_1$ and $K_2$ which illustrates the existence of boundary equilibria when $K_1> x_{1\ell}^* \mbox{ or } K_2> \hat{x}_{2\ell}^*,~\ell=1,2$. The black dots represent real positive $x_{1\ell}^*$ and $\hat{x}_{2\ell}^*$ that satisfy existence of boundary equilibria,   respectively.}
\label{fig1:be}
\end{figure}

Recall from Proposition \eqref{pr1:be} that the interior equilibria $E_{x_1,y_1,y_2}^{l}$  and $E_{y_1,x_2,y_2}^{l}$ $l=1,2$ of Model \eqref{xi0} correspond to the boundary equilibria $E_{1\ell}^b = (x_{1\ell}^*,y_{1\ell}^*,0,y_{2\ell}^*)$ and $E_{2\ell}^b = (0,\hat{y}_{1\ell}^*,\hat{x}_{2\ell}^*,\hat{y}_{2\ell}^*), ~\ell=1,2$ of Model \eqref{2DPatch}. Based on Proposition \eqref{pr1:be} , we could conclude that  Model \eqref{2DPatch}  has four such boundary equilibria. Figures \ref{fig1:be} provide  such an numerical example for the existence of the four boundary equilibria $E_{1\ell}^b = (x_{1\ell}^*,y_{1\ell}^*,0,y_{2\ell}^*)$ and $E_{2\ell}^b = (0,\hat{y}_{1\ell}^*,\hat{x}_{2\ell}^*,\hat{y}_{2\ell}^*)$ under the following parameters:
{\footnotesize
$$s = 0.65,\,\, r_1=1,\,\, r_2=0.54,\,\, d_1=0.45,\,\, d_2=0.105,\,\,K_1=10,\,\,K_2=8,\,\,a_1=0.6,\,\, a_2=0.35,\,\,\rho_1 = 1.75, \,\,\rho_2 = 1.2.$$
}

\begin{figure}[H]
\begin{center}
\subfigure[Effect of dispersal strategy when \newline \hspace*{.45cm} $x_2=0$ and $d_1 = 0.85, a_1=1$.]{\includegraphics[height = 54mm, width = 54mm]{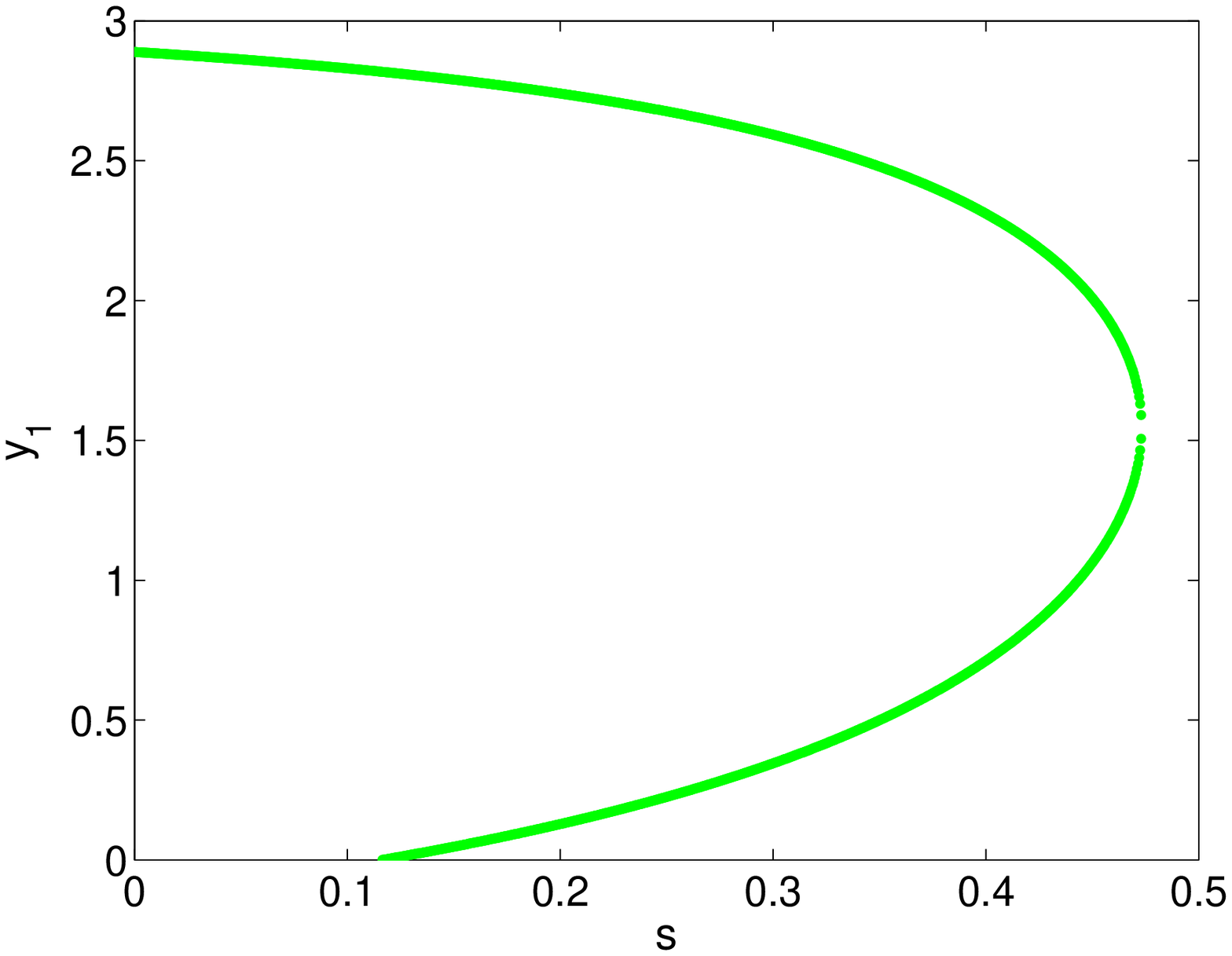}\label{fig_boundaryx20y2_1}}\hspace{.4mm}
\subfigure[Effect of dispersal strategy when \newline \hspace*{.45cm}  $x_1 = 0$ and $d_1 = 0.85, a_1=1$.]{\includegraphics[height = 54mm, width =54mm]{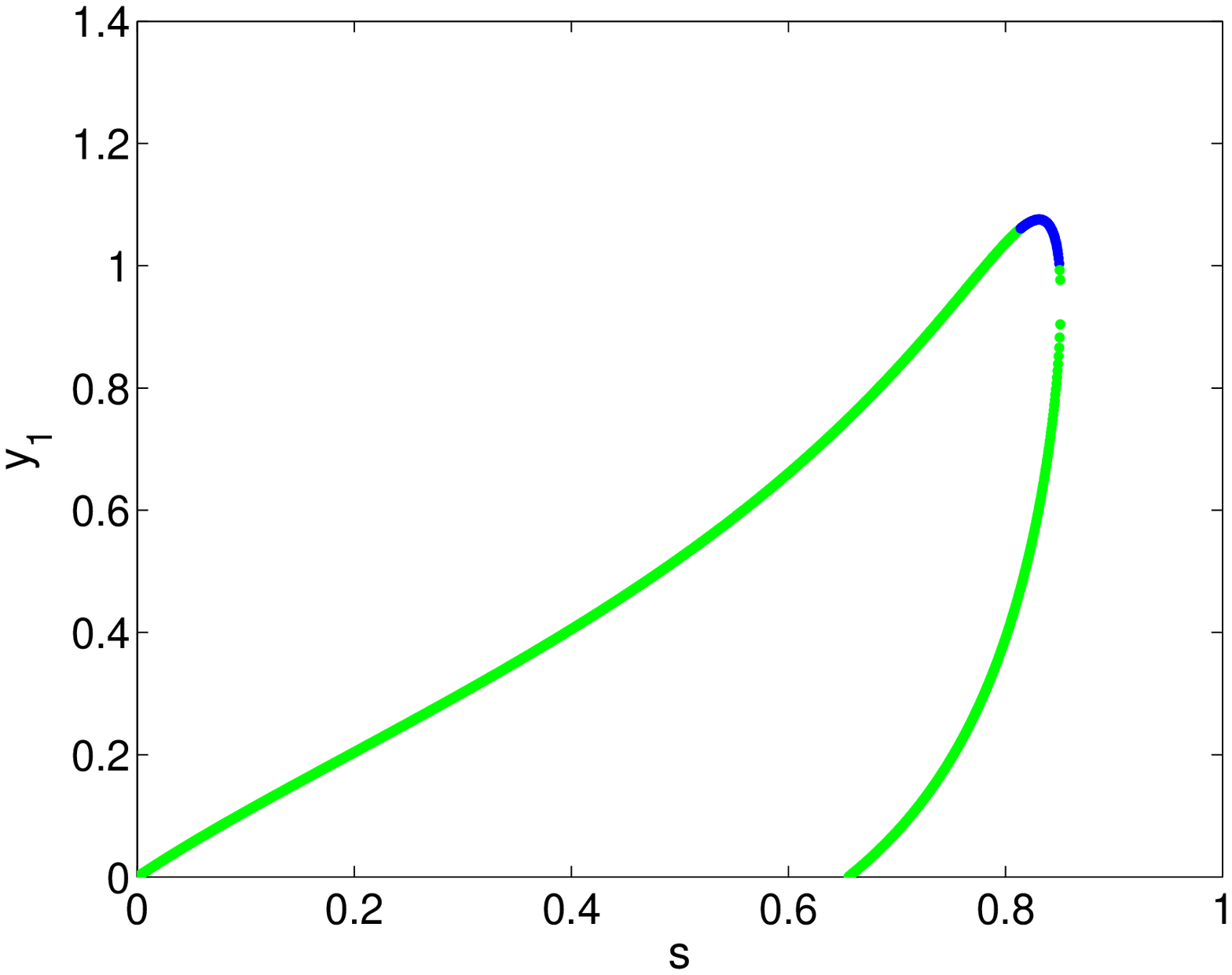}\label{fig_boundaryx10y2_1}}\hspace{.4mm}
\subfigure[Effect of dispersal strategy when \newline \hspace*{.45cm}  $x_1 = 0$ and $d_1 = 2, a_1=2.1$.]{\includegraphics[height = 54mm, width =54mm]{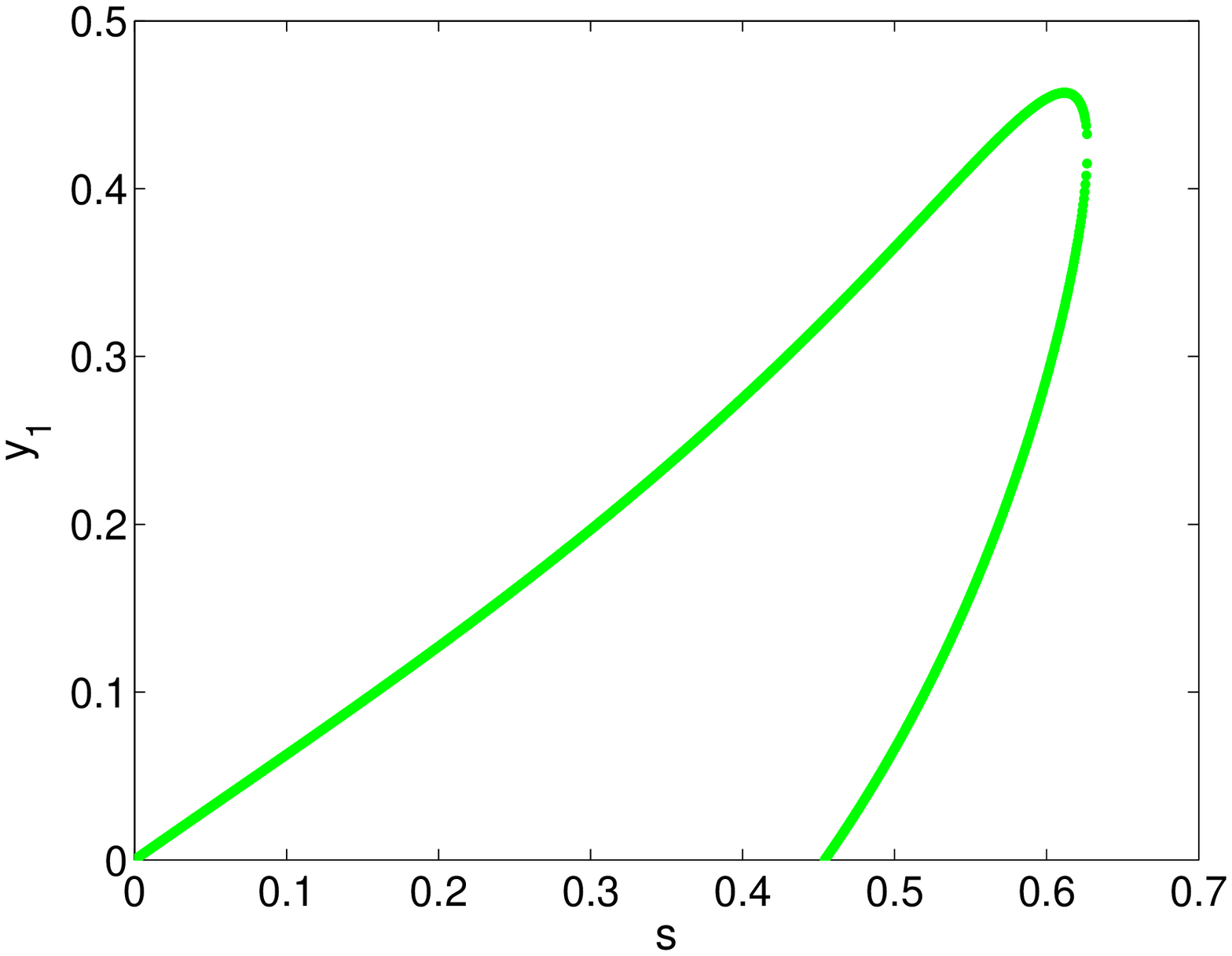}\label{fig_boundaryx10y2_2}}
\end{center}
\vspace{-15pt}
\caption{ One parameter bifurcation diagrams of Model \eqref{2DPatch} with $y$-axis representing the population size of predator at Patch 1 and $x$-axis represent the proportion of predator using the passive dispersal. Figure \ref{fig_boundaryx20y2_1} describes  the number of boundary equilibria $E_{1\ell}^b = (x_{1\ell}^*,y_{1\ell}^*,0,y_{2\ell}^*), \ell=1,2$ from Model \eqref{2DPatch} and their stability with respect to variation in $s$ when $d_1 = 0.85,~a_1=1$. Figure \ref{fig_boundaryx10y2_1} and \ref{fig_boundaryx10y2_2} describes the number of boundary equilibria $E_{2\ell}^b = (0,\hat{y}_{1\ell}^*,\hat{x}_{2\ell}^*,\hat{y}_{2\ell}^*), \ell =1,2$ from Model \eqref{2DPatch} and their change in stability  when $s$ varies from $0$ to $1$ with $d_1 = 0.85,~a_1=1$ and $d_1 = 2,~a_1=2.1$ respectively. Blue represents the sink and green represents the saddle.}
\label{fig3:bb1}
\end{figure}

We continue our study by analyzing the effects of $s$ on the dynamics of the boundary equilibria  $E_{1\ell}^b$ and $E_{2\ell}^b$, $\ell=1,2 $ by adopting the same parameters in generating interior equilibria of Model \eqref{xi0} shown in Figure \ref{fig:intm3_1}, i.e., let $d_1 = 0.85~, a_1 = 1$ and $d_1 = 2,~ a_1 = 2.1$ and 
$$r=1.8,\,\,d_2=0.35,\,\,K_1=10,\,\,K_2=7,\,\,a_2=1.4,\,\,\rho_1 = 1, \,\,\rho_2 = 2.5.$$
Under these parameter values, we have the following two cases that are shown in Figure \ref{fig3:bb1}:
\begin{enumerate}
\item $d_1 = 0.85~, a_1 = 1$: In this case, Model \ref{2DPatch} can have up to three boundary equilibria depending on the values of $s$ (see Figures \ref{fig_boundaryx20y2_1}, \ref{fig_boundaryx10y2_1} and Table \ref{table_Boundary_Model1}).
\item $d_1 = 2,~ a_1 = 2.1$: In this case, Model \ref{2DPatch} can have up to two boundary equilibria depending on the values of $s$ (see Figures \ref{fig_boundaryx10y2_2} and Table \ref{table_Boundary_Model1}).
\end{enumerate}

We recapitulate the following dynamics regarding the effect of $s$ on the equilibria $E_{1\ell}^b$ and $E_{2\ell}^b,~\ell=1,2$: (1) Model \eqref{2DPatch} can have up to four boundary equilibria; (2) These boundary equilibria when exist are locally asymptotically stable or saddle; (3) Large $s$ has a potential to destroy these equilibria. Also, observe the blue line for locally stable and green line for saddle in Figure \ref{fig_InteriorModel3x20_1} as oppose to only green line for saddle in Figure \ref{fig_boundaryx20y2_1}; this results suggest that the additional dimension from the three species Model \eqref{xi0} has a destabilization effect  on the four species Model \eqref{2DPatch}.


{\small\begin{table}[H]
\centering
\begin{tabular}{|c|C{1.4cm}|C{1.4cm}|C{1.4cm}|C{1.4cm}|C{1.4cm}|C{1.4cm}|C{1.4cm}|}
	\hline
\multicolumn{1}{|c|}{} &  \multicolumn{4}{c|}{$\mathbf{a_1= 1 \mbox{  and  } d_1 =0.85 }$ }  & \multicolumn{3}{c|}{$\mathbf{a_1= 2.1 \mbox{  and  } d_1 = 2 }$ }\\
\cline{2-8}	 
{\bf Scenarios }   &  $\mathbf{E_{11}^{b} }$    &   $\mathbf{E_{12}^{b} }$ &$\mathbf{ E_{21}^{b}}$ & $\mathbf{E_{22}^{b} }$ &  $\mathbf{E_{11,12}^{b} }$ &  $\mathbf{ E_{21}^{b}}$& $\mathbf{E_{22}^{b} }$  \\ \hline
$s\leq0.1$ &  Saddle &  \xmark & Saddle &  \xmark &  \xmark & Saddle & \xmark  \\ \hline
$ 0.15 \leq s \leq 0.45$ &   Saddle & Saddle & Saddle &  \xmark &  \xmark & Saddle & \xmark  \\ \hline
$ 0.55 \leq s \leq 0.62$ &  \xmark & \xmark & Saddle &  \xmark &  \xmark & Saddle & Saddle  \\ \hline
$ 0.68 < s <  0.82$ &  \xmark &  \xmark & LAS &  Saddle & \xmark & \xmark& \xmark \\
\hline
$ s \geq 0.82$ & \xmark &  \xmark &  \xmark &  \xmark & \xmark & \xmark & \xmark \\
\hline
\end{tabular}
\caption{{\footnotesize Summary of the effect of the proportion of predators using the passive dispersal on Model \eqref{xi0} From Figures \ref{fig_boundaryx20y2_1}, \ref{fig_boundaryx10y2_1}, and \ref{fig_boundaryx10y2_2}. LAS refers to local asymptotical stability and \xmark \hspace{.01in} implies the equilibrium does not exist.}}
\label{table_Boundary_Model1}
\end{table}}


\subsection{Interior equilibria and stability of Model \eqref{2DPatch} }
Define $p_i(x)=\frac{a_i x}{1+x}$, $q_i(x)=\frac{r_i(K_i-x)(1+x)}{a_iK_i}$, and recall that $\mu_i = \frac{d_i}{a_i-d_i}$. Then from Model \eqref{2DPatch} we have the following equations
\begin{align*}
\frac{dx_i}{dt} &=r_ix_i\left(1-\frac{x_i}{K_i}\right)-\frac{a_i x_iy_i}{(1+x_i)}=\frac{a_i x_i}{1+x_i}\left[\frac{r_i(K_i-x_i)(1+x_i)}{a_iK_i}-y_i\right]=p_i(x_i)\left[q_i(x_i)-y_i\right]. \\
\rho_j\frac{dy_i}{dt} + \rho_i\frac{dy_j}{dt} &= \rho_jy_i[\frac{a_ix_i}{1+x_i}-d_i]+ \rho_iy_j[\frac{a_jx_j}{1+x_j}-d_j] = \rho_jy_i[p_i(x_i)-d_i]+ \rho_iy_j[p_j(x_j)-d_j] 
\end{align*}

Consider $(x_1^*,y_1^*,x_2^*, y_2^*)$ as an interior equilibrium of Model \eqref{2DPatch}, then the following conditions must be satisfied:

\begin{equation}\label{interior-eq1}
\begin{aligned}
& q_i(x_i)-y_i =0 \Leftrightarrow y_i = q_i(x_i) \\[.2cm] 
& \mbox{                     and                 } \\[.2cm]
& \rho_jy_i[p_i(x_i)-d_i]+ \rho_iy_j[p_j(x_j)-d_j] = 0 \Leftrightarrow \rho_jy_i[p_i(x_i)-d_i]= -\rho_iy_j[p_j(x_j)-d_j] 
\end{aligned}
\end{equation}

 which yields the following by substituting the expression of $p_i(x)$ and $q_i(x)$ into \eqref{interior-eq1}
 \begin{equation}\label{interior-eq2}
\begin{aligned}
x_i^2-(\mu_i+K_i)x_i+\underbrace{\mu_iK_i+\frac{a_iK_i}{a_jK_j}\frac{\rho_ir_j}{\rho_jr_i}\frac{(a_j-d_j)}{(a_i-d_i)}(x_j-\mu_j)(x_j-K_j)}_{\phi_i(x_j)} = 0
\end{aligned}
\end{equation}
 
  The equation \eqref{interior-eq2} gives the following nullclines:
 
\begin{equation}\label{interior-eq3}
\begin{aligned}
x_i = \frac{(\mu_i+K_i) \pm \sqrt{(\mu_i+K_i)^2-4\phi_i(x_j)}}{2} = F_i(x_j),~ i,j=1,2,~i\not=j.
\end{aligned}
\end{equation}
The complex form of \eqref{interior-eq3} prevents us to obtain the explicit solutions of the interior equilibria of Model \eqref{2DPatch}. We are going to explore the symmetric interior equilibrium for the symmetric Model \eqref{2DPatch} where we say that Model \eqref{2DPatch} is symmetric if $~ a_1 = a_2 =a,~d_1=d_2=d,~K_1=K_2=K,~r_1=r_2=r.$ Now we have the following theorem:\\



\begin{figure}[H]
\begin{center}
\subfigure[$s$ V.S. $y_1$ for the effect of $s$ when Model \eqref{2DPatch} \newline \hspace{1cm} is symmetric with $\rho_1=1.72$ and $\rho_2 = 13$]{\includegraphics[height = 55mm, width = 65mm]{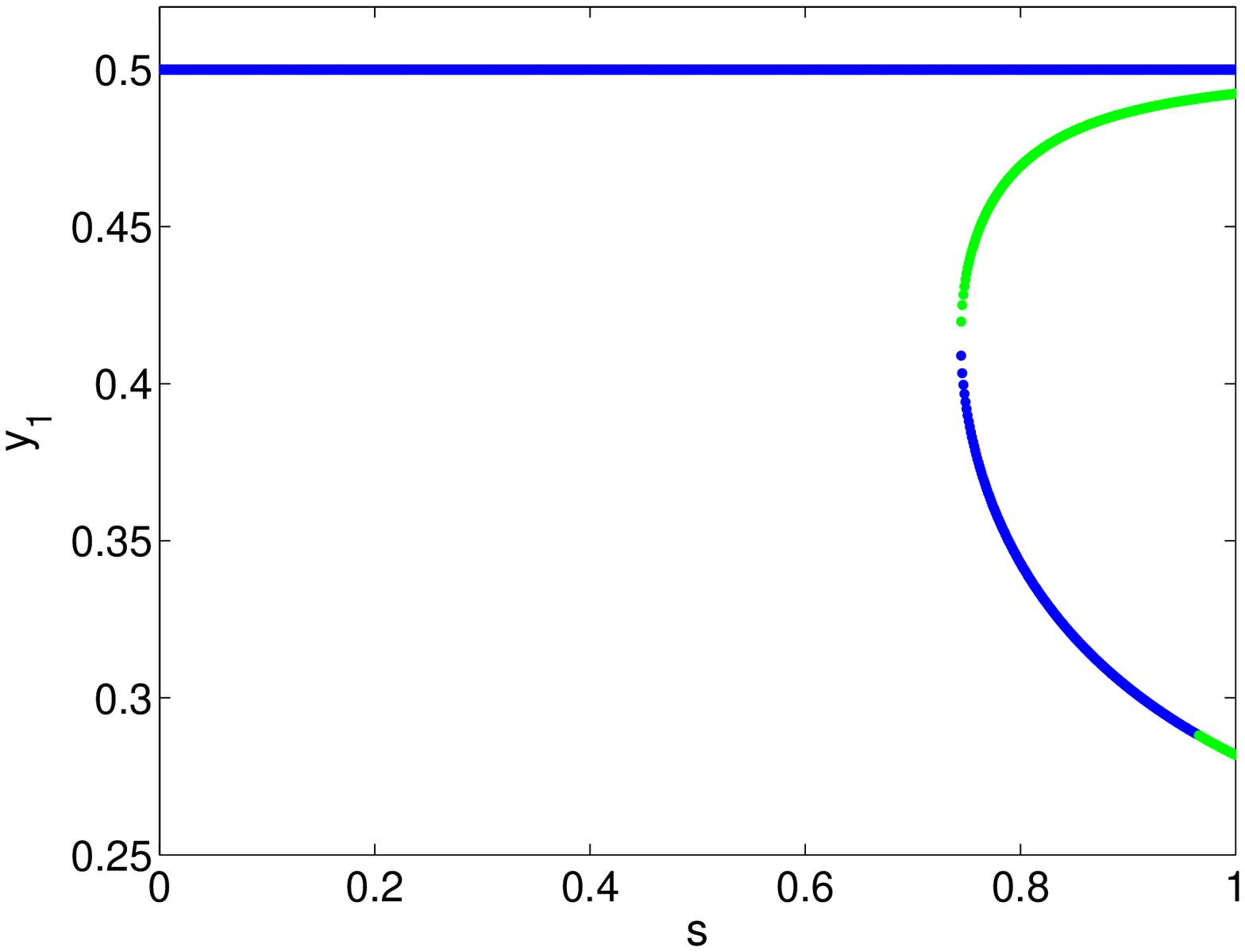}\label{fig_InteriorSyModel1_Sy1}}\hspace{2mm}
\subfigure[$s$ V.S. $\rho_1$ for the number of interior equilibria when Model \eqref{2DPatch} is symmetric and $\rho_2 = 13$]{\includegraphics[height = 55mm, width = 65mm]{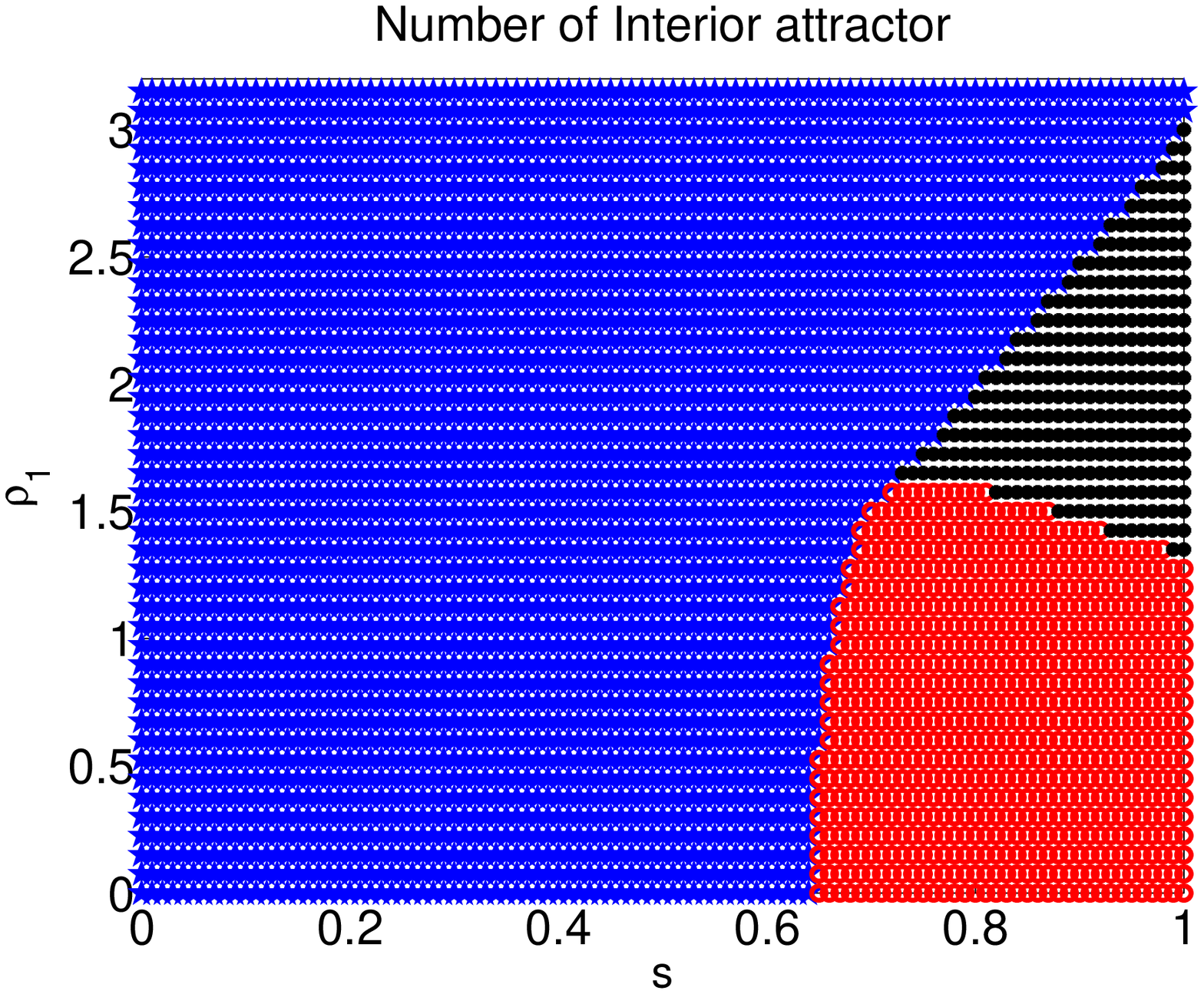}\label{fig_InteriorSyModel1_rho1S}}
\end{center}
\vspace{-15pt}
\caption{One and two parameter bifurcation diagrams of symmetric Model \eqref{2DPatch} with $y$-axis representing the population size of predator at Patch 1 in figure \ref{fig_InteriorSyModel1_Sy1}. We used the following parameters $r=1$, $d=5$, $K=10$, and $a=6$. Figure \ref{fig_InteriorSyModel1_Sy1} describes the number of interior equilibria and their change in stability  when $s$ varies from $0$ to $1$. Blue line represents sink and green line represents saddle in Figure \ref{fig_InteriorSyModel1_Sy1}. Figure \ref{fig_InteriorSyModel1_rho1S} describes how the number of interior equilibria change for different values of dispersal strategy $s$ and dispersal rate $\rho_1$. Black region have three interior equilibria; red regions have two interior equilibria; and blue regions have one interior equilibrium in Figure \ref{fig_InteriorSyModel1_rho1S}.}
\label{fig:M1IntSy}
\end{figure}

\begin{theorem}\label{th34:interior}[The symmetric interior equilibrium and the stability]
Suppose that Model \eqref{2DPatch} is symmetric with $r=1$. We denote
 $$\mu=\frac{d}{a-d},\mbox{   and  }  \nu=\frac{(K-\mu)(1+\mu)}{aK}.$$
 Then $E=(\mu,\nu,\mu,\nu)$ is an unique symmetric interior equilibrium for Model \eqref{2DPatch}. Moreover, $E$ is locally asymptotically stable if $\frac{K-1}{2}<\mu<K$ while it is unstable if $\mu<\frac{K-1}{2}$ for $s \in [0,1]$.
 \end{theorem}
\noindent\textbf{Notes:} Theorem \eqref{th34:interior} implies the symmetric Model \eqref{2DPatch} has an unique symmetric interior equilibrium of the form $E=(\mu,\nu,\mu,\nu)$. The related results imply that dispersal of predators and $s$ has no effect on the local stability of this symmetric interior equilibria when it exist since $\frac{K-1}{2}<\mu<K$ does not depend on $\rho_i,i=1,2$ or $s$. We note that Model \eqref{2DPatch} can have two additional interior equilibria in the symmetric case which can be locally stable or saddle depending on the value of $s$  (see green line for saddle and blue line for locally stable in Figures \ref{fig_InteriorSyModel1_Sy1} which correspond to the additional two boundary equilibria of Model \eqref{2DPatch} in the symmetric case). We consider the following fixed symmetric parameters:


$$r_1=r_2=r=1,\,\,d_1=d_2=d=5,\,\,K_1=K_2=K=10,\,\,a_1=a_2=a=6.$$

According to the bifurcation diagrams in Figures \ref{fig_InteriorSyModel1_Sy1} and \ref{fig_InteriorSyModel1_rho1S}, Model \eqref{2DPatch} can have up to three interior equilibria in the symmetric case. It seems that the larger value of $s$ can create two additional asymmetric interior equilibria which can be saddle or locally stable, thus generate bistability between two different interior attractors (See blue lines in Figure \ref{fig_InteriorSyModel1_Sy1} when $0.78\leq s \leq 0.92$). The local stability of $E=(\mu,\nu,\mu,\nu)$ does not depend on $s$ as illustrated in Theorem \ref{th34:interior}. \\

\noindent\textbf{Summary:} In addition to the summary of our analysis listed in Table \eqref{table_stability}, we summarize the following dynamics of Model \eqref{2DPatch} base on mathematical analysis and bifurcation diagrams from our study:
\begin{enumerate}
\item The four basic boundary equilibria $E_{0000}$, $E_{K_1000}$, $E_{00K_20,}$, $E_{K_10K_20}$ always exist where $E_{0000}$, $E_{K_1000}$, $E_{00K_20,}$ are always saddle while $E_{K_10K_20}$  is locally asymptotically stable if the two inequalities \ref{Ek10k20st} are satisfied. Large dispersal of predators can stabilize the boundary equilibrium $E_{K_10K_20}$ when $s \in (0,1]$. However, the value of $s$ has no effects on the global stability of the boundary equilibrium $E_{K_10K_20}$.
\item Model \eqref{2DPatch} can have up to four other boundary equilibria $E_{1\ell}^b = (x_{1\ell}^*,y_{1\ell}^*,0,y_{2\ell}^*)$ and $E_{2\ell}^b = (0,\hat{y}_{1\ell}^*,\hat{x}_{2\ell}^*,\hat{y}_{2\ell}^*)$ for $\ell = 1,2$. The number of these boundary equilibria and the stability could be affected by the dispersal strength $\rho_i,i=1,2$ and the values of $s$. For example, the large values of $s$ can destroy these boundary equilibria.
\item In the symmetric case, Model \eqref{2DPatch} may potentially have three interior equilibria including $E = (\mu,\nu,\mu,\nu)$ from Theorem \ref{th34:interior} when $a > d$. Although the local stability of $E$ does not depend on $s$, the large value of $s$ can generate the two additional asymmetric interior equilibria, hence, create multiple interior attractors.\\
\end{enumerate}

Define $\mu_i=\frac{d_i}{a_i-d_i},\nu_1=\frac{(K_1-\mu_1)(1+\mu_1)}{a_1K_1}, \nu_2=\frac{r(K_2-\mu_2)(1+\mu_2)}{a_2K_2}$, $\hat{\mu_i}=\frac{\hat{d_i}}{a_i-\hat{d_i}}$, $ \hat{\nu_i}=q_i(\hat{\mu_i})=\frac{r_i(K_i-\hat{\mu_i})(1+\hat{\mu_i})}{a_iK_i}$, $\hat{\nu}_j^i=\frac{\rho_j\hat{\nu_i}}{d_j+\rho_j}$ where $\hat{d_i}=d_i +\frac{\rho_id_j}{d_j+\rho_j}$ $i,j=1,2,\,i\neq j$ and $E_{12}^{b*}=E_{\mu_1\nu_1K_20},\,E_{22}^{b*}=E_{K_10\mu_2\nu_2}$. Then the boundary dynamics for $s=0,1$ from the work of \citep{jansen2001dynamics, kang2014dispersal} and $s\in(0,1)$ from our current work is summarize in Table \ref{table_stability}.\\


{\small\begin{table}[H]
\centering
\begin{tabular}{|c|p{3.5cm}|p{4.2cm}|p{4.2cm}|}\hline
\multicolumn{1}{|c|}{} &  \multicolumn{3}{c|}{{\bf Existence condition, Local and Global stability of Model \eqref{2DPatch}} } \\
\cline{2-4}
{\bf Scenarios}   &   $\mathbf{ s = 0}$    &   $\mathbf{ s \in (0,1)}$ & $\mathbf{ s = 1}$    \\
\hline
$E_{0000},$ $E_{K_1000},$ $E_{00K_20}$ & Always exist and always saddle & Always exist and always saddle &  Always exist and always saddle\\
\hline
$E_{K_10K_20}$  & Always exist; LAS and GAS if $\mu_i>K_i$ for both $i=1,2$& Always exist; GAS if $\mu_i>K_i$ for both $i=1,2$; while LAS if Equations \ref{Ek10k20st} are satisfied & Always exist; GAS if $\mu_i>K_i$ for both $i=1,2$; LAS if condition condition (1) is satisfied \\
\hline
\pbox{20cm}{$E_{1\ell}^b$ ($x_i=0$), $\ell = 1,~2$\\ $i =1,2$} & Do not exist & One or two exist if {\footnotesize $\frac{3\beta_j}{\mu_j+K_j}<\alpha_j<(\mu_j+K_j)^2$ } with {\footnotesize $i,j=1,2$, $i\not=j$}; Can be locally asymptotically stable or saddle as shown in Figures \ref{fig_boundaryx20y2_1}, \ref{fig_boundaryx10y2_1}, \ref{fig_boundaryx10y2_2} & Exist if {\footnotesize $0<\hat{\mu}_i<K_i$}; LAS if  {\footnotesize $\frac{K_i-1}{2}<\widehat{\mu}_i<K_i$} and {\footnotesize $r_j<a_j\hat{\nu}_j^i$}. GAS if {\footnotesize $\frac{K_i-1}{2}<\widehat{\mu}_i<K_i$} and  {\footnotesize $\frac{r_j(K_j+1)^2}{4a_jK_j}<\widehat{\nu}_i^j$}, {\footnotesize $i,j=1,2$, $i\not=j$}.\\
\hline

$E_{i2}^{b*}$, $i,j=1,2,\,i\neq j$ & Exist if {\footnotesize $0<\mu_i<K_i$}; LAS if $\frac{K_i-1}{2}<\mu_i<K_i$ and condition (2) is satisfied & Do not exist  & Do not exist \\
\hline

\multicolumn{1}{|c|}{{\bf Condition 1:}} &  \multicolumn{3}{c|}{ {\tiny $\mathlarger{\sum}_{i=1}^{2}\left[ \frac{(a_i-d_i)(\mu_i-K_i)}{1+K_i} + \rho_i \right] > 0$} and {\tiny $\left[\frac{(a_1-d_1)(\mu_1-K_1)}{1+K_1}\right]\left[\rho_2+\frac{(a_2-d_2)(\mu_2-K_2)}{1+K_2}\right]+\rho_1\left[\frac{(a_2-d_2)(\mu_2-K_2)}{1+K_2}\right] > 0$ } } \\
\multicolumn{1}{|c|}{{\bf Condition 2:}} &  \multicolumn{3}{c|}{ {\tiny $0<\frac{d_i}{a_j-d_i}<K_j<\mu_j$} and {\tiny $\rho_j<\frac{d_j-K_j(a_j-d_j)}{\nu_i\left[K_j(a_j-d_i)-d_i\right]}$; $i,j=1,2$, $i\not=j$} } \\
\hline
\end{tabular}
\caption{Summary of the local and global dynamic of Model \eqref{2DPatch}. LAS refers to the local asymptotical stability, and GAS refers to the global stability.}
\label{table_stability}
\end{table}}


\section{Effects of dispersal strategies on the prey-predator population dynamics}\label{sec_effects_of_dispersal_s}
In order to get more insights into the dynamics of Model \eqref{2DPatch}, we perform bifurcation analysis in this section. We fixed the following parameters for most of the simulations 
$$r=1.8,\,\,d_2=0.35,\,\,K_1=10,\,\,K_2=7,\,\,a_2=1.4,\,\,\rho_1 = 1, \,\,\rho_2 = 2.5.$$
 and consider these two cases: $d_1 = 0.85,~a_1=1$ and $d_1 = 2,~a_1=2.1$. According to the dynamics of the subsystem Model \eqref{Onepatch} provided in Section \ref{ModelDerivation}, we know that in the absence of dispersal, Patch 1 has global stability at $(10,0)$ if $\frac{d_1}{a_1-d_1} > 10$ (e.g., when $d_1 = 2,~a_1=2.1$) and it has global stability at its unique interior $\left(\frac{d_1}{a_1-d_1},\frac{\left(10-\frac{d_1}{a_1-d_1}\right)\left(1+\frac{d_1}{a_1-d_1}\right)}{10a_1}\right)$ if $4.5 < \frac{d_1}{a_1-d_1} < 10$ (e.g., when $d_1 = 0.85,~a_1=1$); while Patch 2 has a unique stable limit cycle since $d_2=0.35,,\,\,K_2=7,\,\,a_2=1.4$. \\
 


We implement one and two parameters bifurcation diagrams to obtain  insights into the dynamical patterns of the asymmetric two patch Model \eqref{2DPatch} in the following way:

\begin{enumerate}
\item  $d_1 = 0.85$ and $a_1 = 1$:  In the absence of dispersal, the uncoupled two patch model is unstable at the interior equilibrium $(5.67,288.89, 0.33,80)$. However, in the presence of the dispersal, Figure \ref{fig_InteriorModel1_Sy1_1} (blue regions) suggest that the intermediate values of $s$ can stabilize the dynamics while the large values of $s$ with certain dispersal strengths could generate multiple interior equilibria (up to three interior equilibria), thus lead to multiple attractors potentially. Moreover, two dimensional bifurcation diagram shown in Figure \ref{fig_InteriorModel1_rho1S_1}
 suggest that the large values of $s$ combined with the small or large dispersal strength $\rho_1$ in Patch 1 can destroy the interior equilibria (see white regions in Figure \ref{fig_InteriorModel1_rho1S_1}) with consequences that prey in one patch may go extinct but predator persists in each patch. Table \eqref{table_Interior_Model1} provides a more details description on the existence and stability of the interior equilibria of Model \eqref{2DPatch}. 

\begin{figure}[H]
\begin{center}
\subfigure[$s$ V.S. $y_1$ for the effect of $s$ when $d_1 = 0.85$, $a_1 = 1$, $\rho_1 = 1$, and $\rho_2 = 2.5$]{\includegraphics[height = 50.mm, width = 65mm]{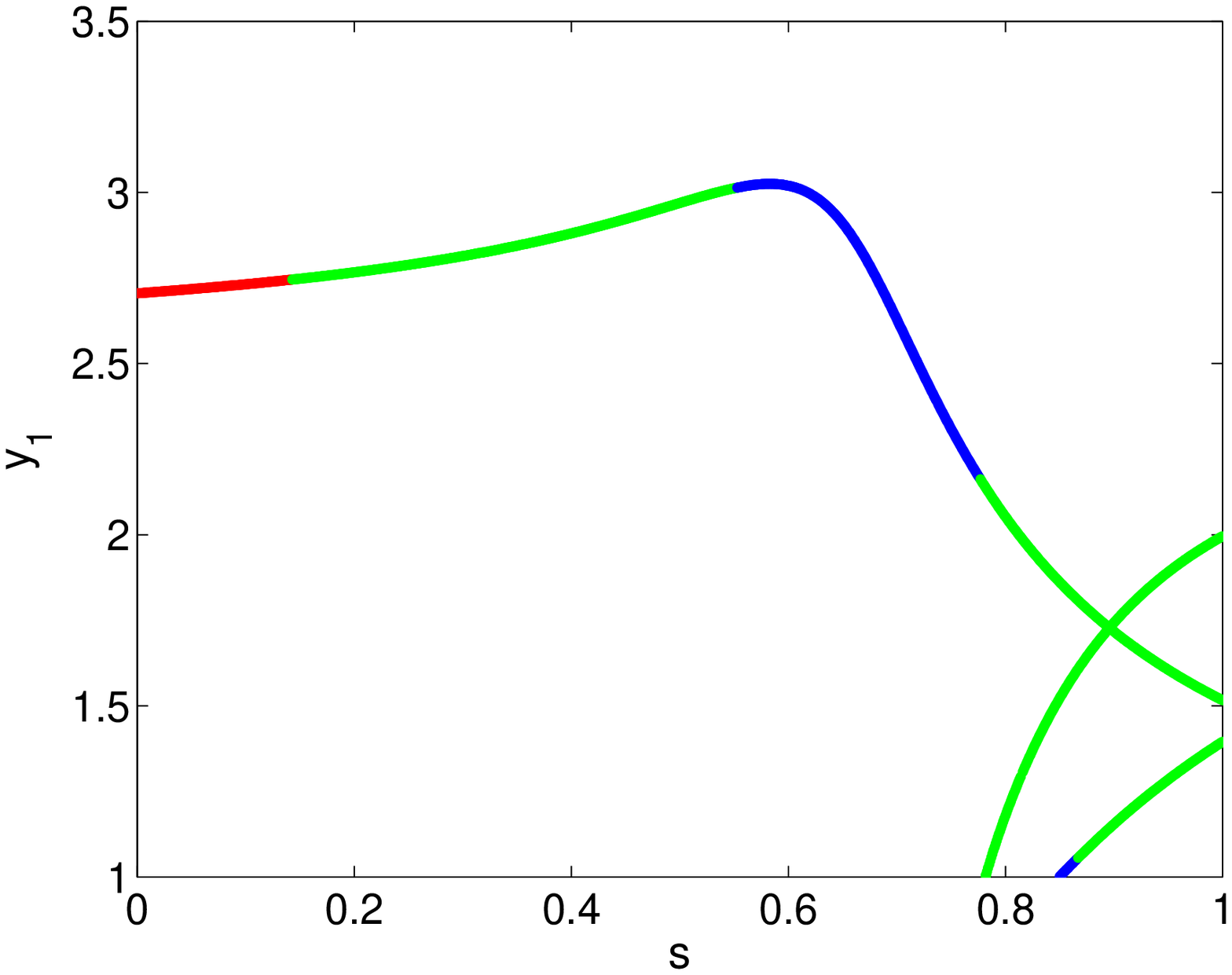}\label{fig_InteriorModel1_Sy1_1}}\hspace{2mm}
\subfigure[$s$ V.S. $y_1$ for the number of interior equilibria when $d_1 = 0.85$, $a_1 = 1$, and $\rho_2 = 2.5$]{\includegraphics[height = 50.mm, width = 65mm]{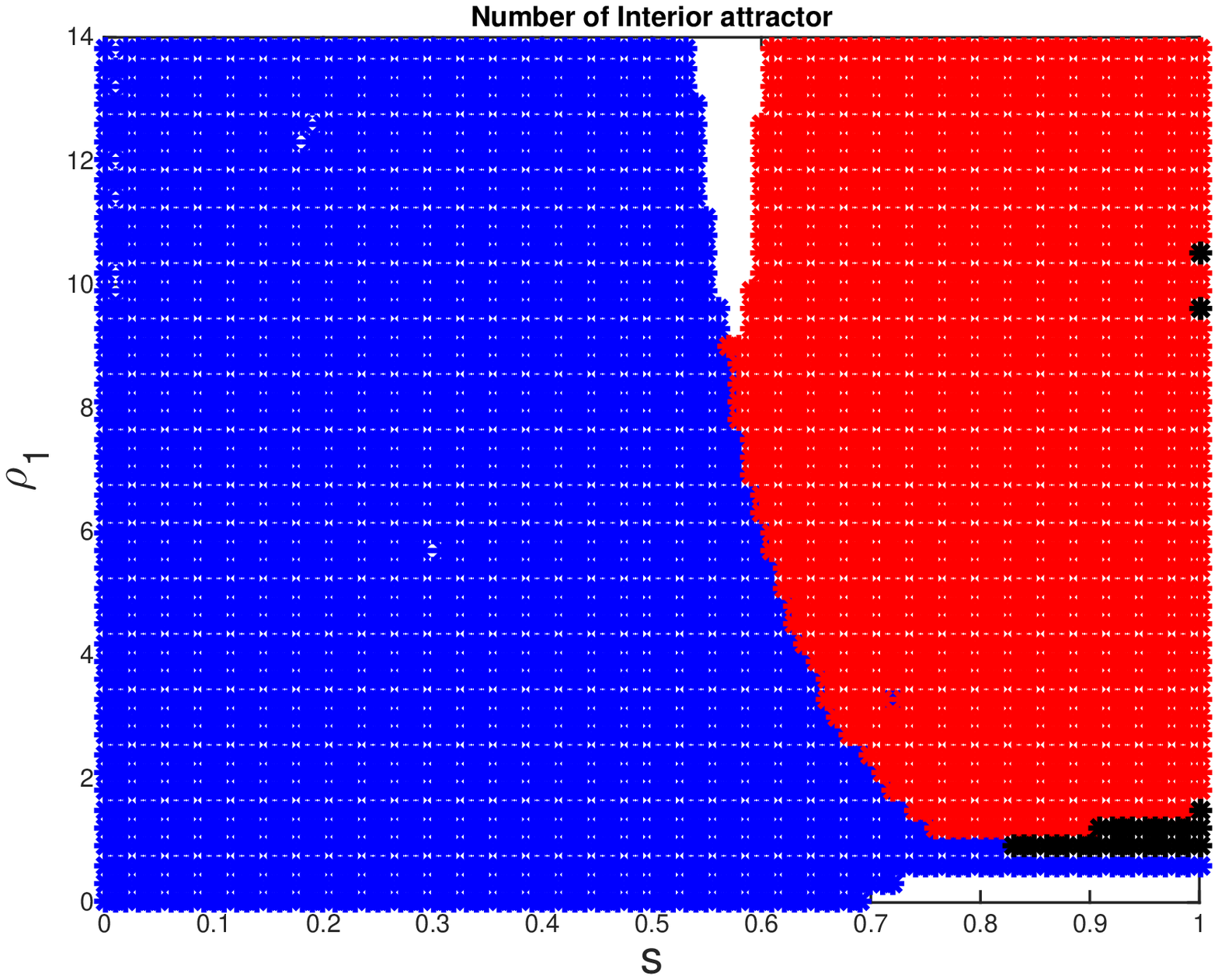}\label{fig_InteriorModel1_rho1S_1}}
\end{center}
\vspace{-15pt}
\caption{ One and two parameter bifurcation diagrams of Model \eqref{2DPatch} with $y$-axis representing the population size of predator at Patch 1 in Figure \ref{fig_InteriorModel1_Sy1_1}. The following parameters are used: $r=1.8$, $d_2=0.35$, $K_1=10$, $K_2=7$, and $a_2=1.4$. Figure \ref{fig_InteriorModel1_Sy1_1} describes the number of interior equilibria and their change in stability  when $s$ varies from $0$ to $1$. Blue line represents sink, green line represents saddle, and red line represents source in Figure \ref{fig_InteriorModel1_Sy1_1}. Figure \ref{fig_InteriorModel1_rho1S_1} describes how the number of interior equilibria change for different values of $s$ and dispersal rate $\rho_1$. Black region have three interior equilibria; red regions have two interior equilibria; blue regions have one interior equilibrium, and white regions have no interior equilibria in Figure \ref{fig_InteriorModel1_rho1S_1} .}
\label{fig1:M1int1}
\end{figure}

{\small\begin{table}[H]
\centering
\begin{tabular}{|c|C{1.75cm}|C{1.75cm}|C{1.75cm}|C{1.75cm}|C{1.75cm}|C{1.75cm}|C{1.75cm}|}
	\hline
\multicolumn{1}{|c|}{} &  \multicolumn{3}{c|}{$\mathbf{a_1= 1 \mbox{  and  } d_1 =0.85 }$ }  & \multicolumn{3}{c|}{$\mathbf{a_1= 2.1 \mbox{  and  } d_1 = 2 }$ }\\
\cline{2-7}	 
{\bf Scenarios }   &  $\mathbf{E_{x_1y_1x_2y_2}^{1} }$    &   $\mathbf{E_{x_1y_1x_2y_2}^{2} }$ &$\mathbf{ E_{x_1y_1x_2y_2}^{3}}$ & $\mathbf{E_{x_1y_1x_2y_2}^{1} }$ &  $\mathbf{E_{x_1y_1x_2y_2}^{2} }$ &  $\mathbf{ E_{x_1y_1x_2y_2}^{3}}$
  \\ \hline
$s\leq0.07$ &  Source &  \xmark & \xmark &  Saddle &  Source & LAS   \\ \hline
$ 0.9 \leq s \leq 0.15$ &   Source & \xmark & \xmark &  Saddle &  Saddle & LAS  \\ \hline
$ 0.2 \leq s \leq 0.43$ &  Saddle & \xmark & \xmark &  LAS &  Saddle & Saddle  \\ \hline
$ 0.55 \leq s \leq  0.68$ &  LAS &  \xmark &  \xmark &  LAS & \xmark & \xmark \\
\hline
$ 0.78 \leq s \leq  0.82$ &  Saddle &  Saddle & \xmark &  Saddle & \xmark & \xmark \\
\hline
$ 0.83 \leq s \leq  0.84$ &  Saddle &  Saddle & LAS &  Saddle & \xmark & \xmark \\
\hline
$ s \geq 0.84$ & Saddle &  Saddle &  Saddle &  Saddle & \xmark & \xmark  \\
\hline
\end{tabular}
\caption{{\footnotesize Summary of the effect of the proportion of predators using the passive dispersal on the interior equilibria of Model \eqref{2DPatch} From Figures \ref{fig_InteriorModel1_Sy1_1}, and \ref{fig_InteriorModel1_Sy1_2}. LAS refers to local asymptotical stability, \xmark \hspace{.01in} implies the equilibrium does not exist, and $E_{x_1y_1x_2y_2}^{i},i=1,2,3$ are the three possible interior equilibria of Model \eqref{2DPatch}.}}
\label{table_Interior_Model1}
\end{table}}
 

\item $d_1 = 2$ and $a_1=2.1$: In the absence of dispersal, the uncoupled two patch model has extinction of predator in Patch 1 and is unstable at the boundary equilibrium $(10,0, 0.33,80)$. However, in the presence of the dispersal, Figure \ref{fig_InteriorModel1_Sy1_2} (blue regions) suggest that the intermediate values of $s$ can stabilize the dynamics while the small values of $s$ with certain dispersal strengths could generate multiple interior equilibria (up to three interior equilibria), thus lead to multiple attractors potentially. Moreover, two dimensional bifurcation diagram shown in Figure \ref{fig_InteriorModel1_rho1S_2} suggest that the large values of $s$ combined with the large dispersal strength $\rho_1$ in Patch 1 can destroy the interior equilibria (see white regions in Figure \ref{fig_InteriorModel1_rho1S_1}) with consequences that prey in one patch may go extinct but predator persists in each patch. A more detail dynamic from Figure \ref{fig_InteriorModel1_rho1S_2} is presented in Table \eqref{table_Interior_Model1}. 


\begin{figure}[H]
\begin{center}
\subfigure[$s$ V.S. $y_1$ for the effect of $s$ when $d_1 = 2$, $a_1 = 2.1$, $\rho_1 = 1$, and $\rho_2 = 2.5$]{\includegraphics[height = 55mm, width = 65mm]{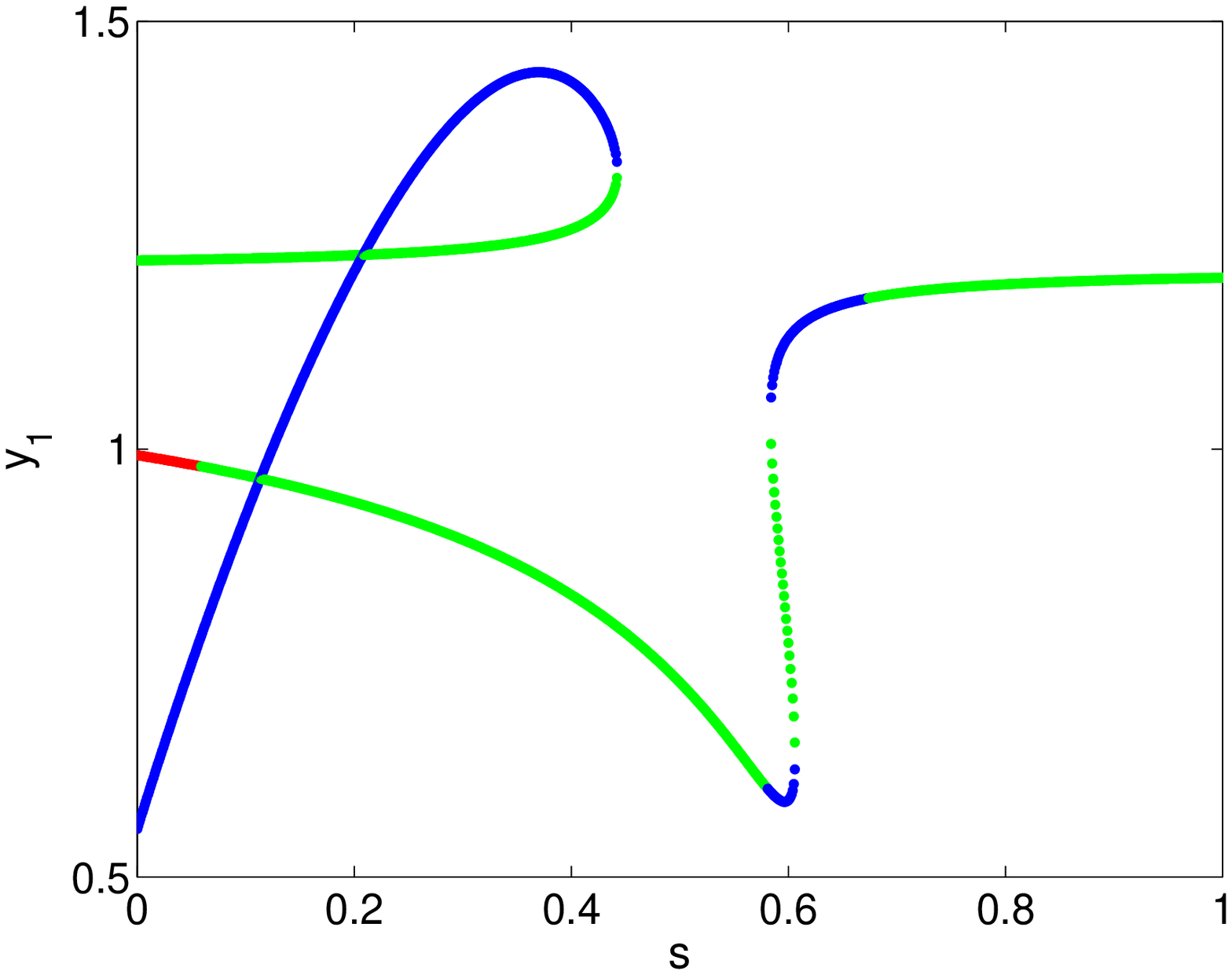}\label{fig_InteriorModel1_Sy1_2}}\hspace{2mm}
\subfigure[$s$ V.S. $\rho_1$ for the number of interior equilibria when $d_1 = 2$, $a_1 = 2.1$,  and $\rho_2 = 2.5$]{\includegraphics[height = 55mm, width = 65mm]{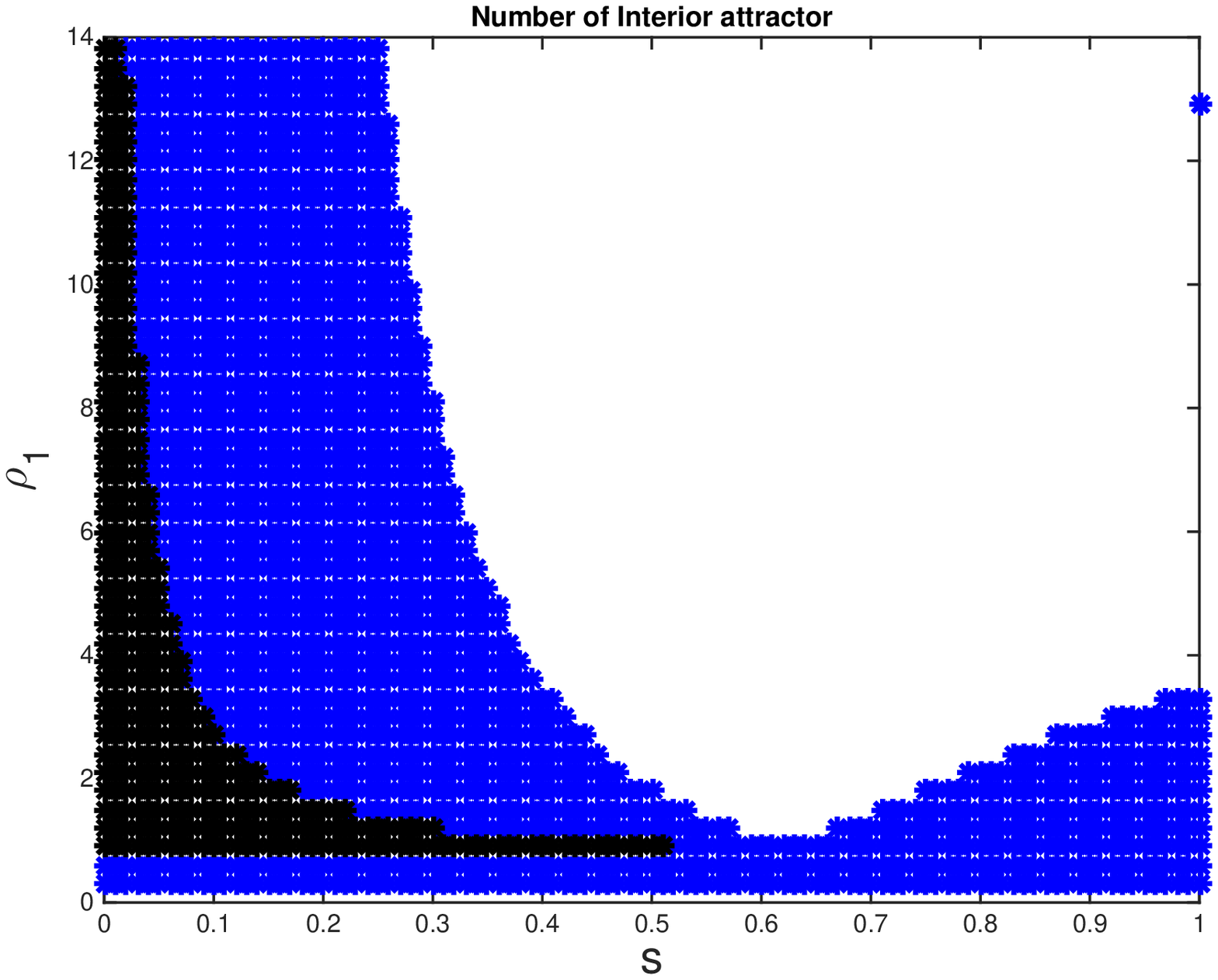}\label{fig_InteriorModel1_rho1S_2}}
\end{center}
\vspace{-15pt}
\caption{One and two parameter bifurcation diagrams of Model \eqref{2DPatch} with $y$-axis representing the population size of predator at Patch 1 in figure \ref{fig_InteriorModel1_Sy1_2}. We used the following parameters $r=1.8$, $d_2=0.35$, $K_1=10$, $K_2=7$, and $a_2=1.4$. Figure \ref{fig_InteriorModel1_Sy1_2} describes the number of interior equilibria and their change in stability  when $s$ varies from $0$ to $1$ under the parameters $d_1 =2$ and $a_1=2.1$. Blue represents sink, green represents saddle, and red represents source. Figure \ref{fig_InteriorModel1_rho1S_2} describes how the number of interior equilibria change for different values of dispersal strategy $s$ and dispersal rate $\rho_1$. Black region have three interior equilibria; red regions have two interior equilibria; blue regions have one interior equilibrium,  and white regions have no interior equilibria.}
\label{fig2:M1Int2}
\end{figure}


\item Two parameter bifurcation diagrams of the relative dispersal rate $\rho_2$ versus the dispersal strategy $s$ for both scenarios of $d_1 = 0.85,~a_1 = 1$ (Figure \ref{fig_InteriorModel1_rho2S_1}) and $d_1 = 2,~a_1=2.1$ (Figure \ref{fig_InteriorModel1_rho2S_2}). For both cases, the large $s$ combined with the large dispersal strength in Patch 2, i.e., $\rho_2$, can destroy the interior equilibrium (see white regions in Figures \ref{fig_InteriorModel1_rho2S_1} and \ref{fig_InteriorModel1_rho2S_2} for $s > 0.6$); while the small $s$ (for $d_1 = 0.85,~a_1 = 1$) and the large value of $s$ (for $d_1 = 2,~a_1=2.1$) could generate multiple interior equilibria (see black region for three interior equilibria and red region for two interior equilibria in Figure \ref{fig_InteriorModel1_rho2S_1} and  \ref{fig_InteriorModel1_rho2S_2}).\\


\begin{figure}[H]
\begin{center}
\subfigure[$s$ V.S. $\rho_2$ for the number of interior equilibria when $d_1 = 0.85$ and $a_1 = 1$]{\includegraphics[height = 50mm, width =60mm]{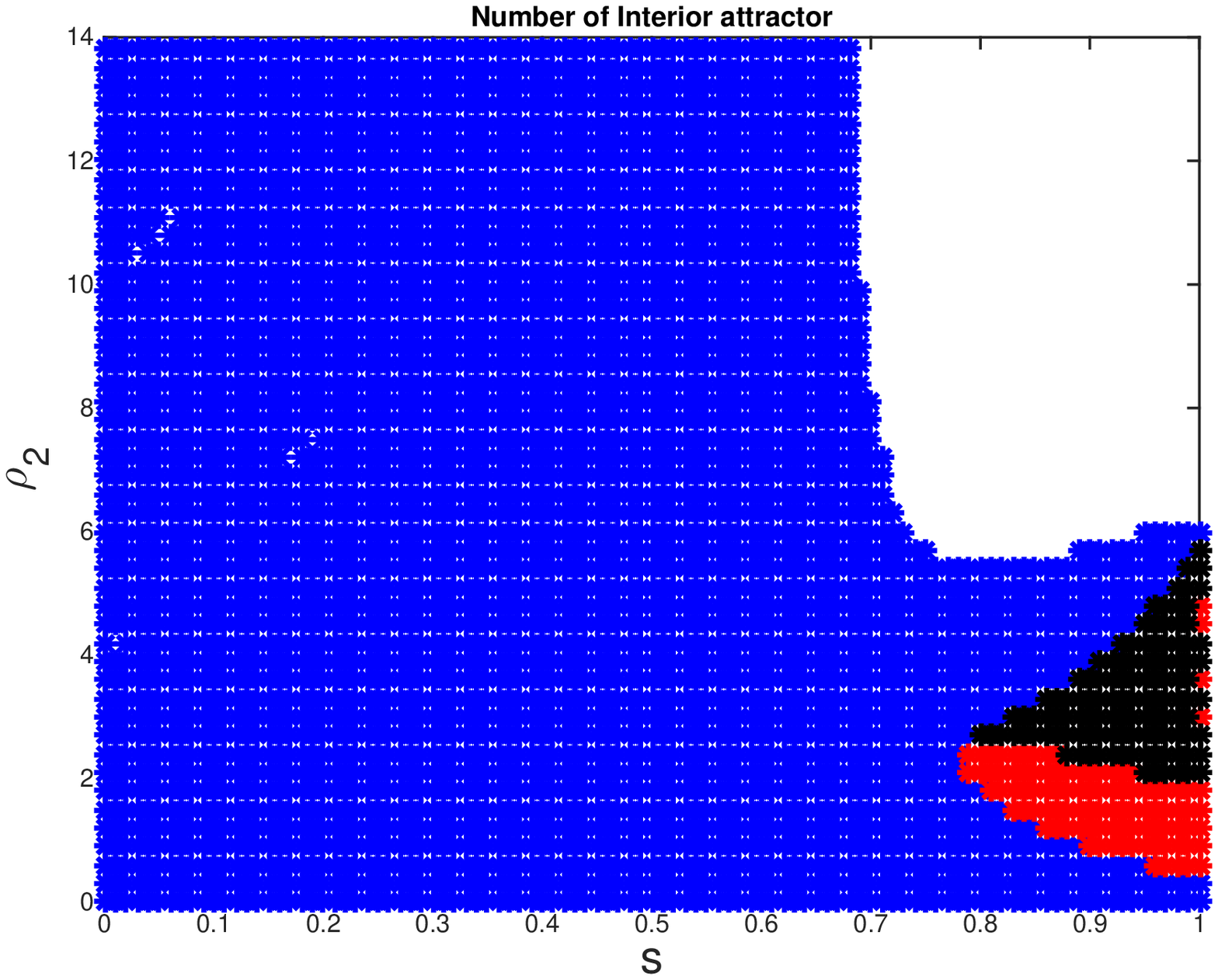}\label{fig_InteriorModel1_rho2S_1}}\hspace{2mm}
\subfigure[$s$ V.S. $\rho_2$ for the number of interior equilibria when $d_1 = 2$ and $a_1 = 2.1$]{\includegraphics[height = 50mm, width = 60mm]{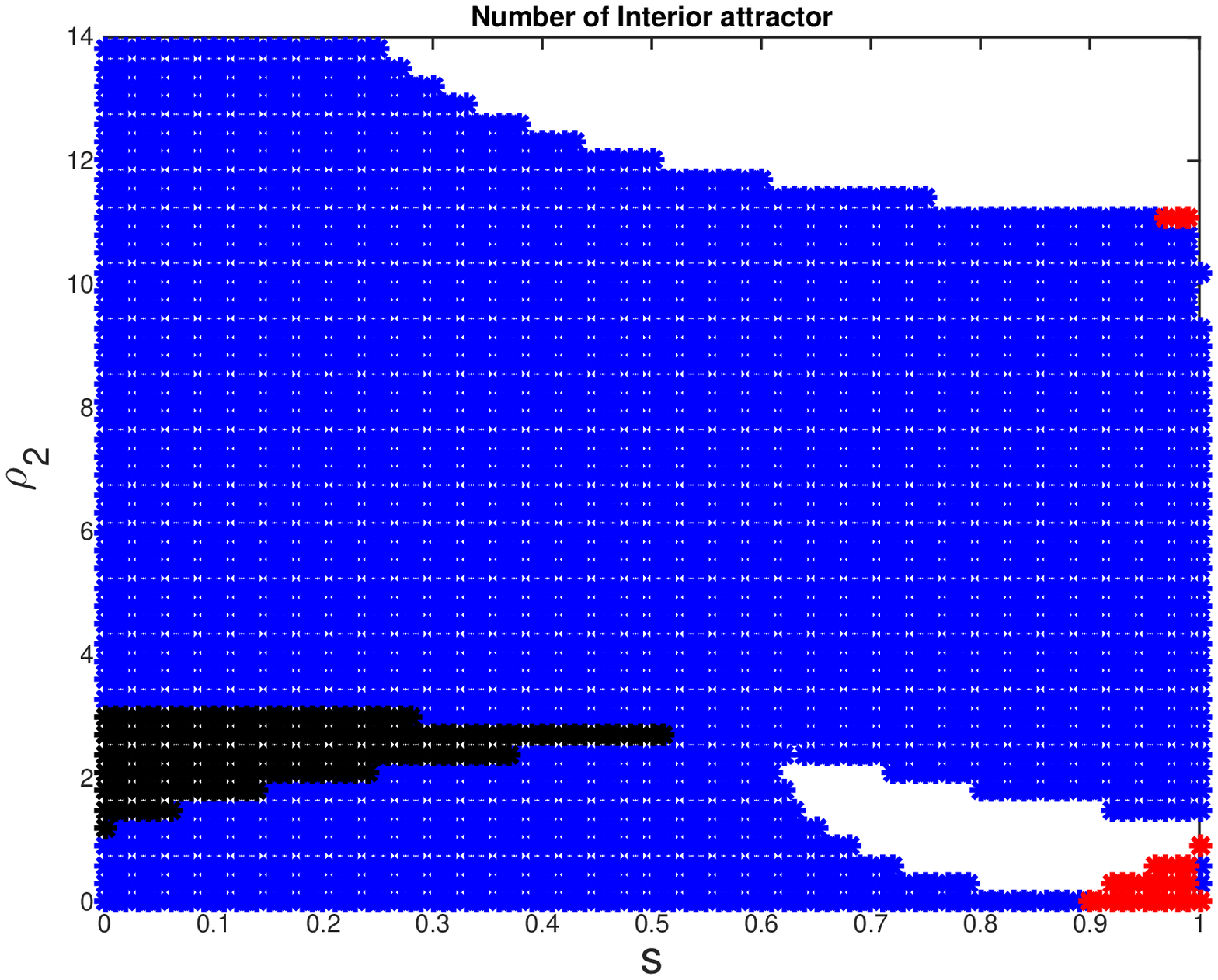}\label{fig_InteriorModel1_rho2S_2}}
\end{center}
\vspace{-15pt}
\caption{Two parameters bifurcation diagrams of Model \eqref{2DPatch} with $y$-axis representing the relative dispersal rate $\rho_2$ and $x$-axis represent the strength of dispersal mode $s$. The following parameters are used: $r=1.8$, $d_2=0.35$, $K_1=10$, $K_2=7$, $a_2=1.4$, and $\rho_1 = 1$. Both figures \ref{fig_InteriorModel1_rho2S_1} and \ref{fig_InteriorModel1_rho2S_2} describes how the number of interior equilibria change for different values of dispersal strategy $s$ and dispersal rate $\rho_2$ where the parameters $d_1 =0.85,~a_1=1$ are used for the left figure \ref{fig_InteriorModel1_rho2S_1} while $d_1 =2,~a_1=2.1$ are used for the right figure \ref{fig_InteriorModel1_rho2S_2} in addition to the fixed parameters. Black region have three interior equilibria; red regions have two interior equilibria; blue regions have one interior equilibrium, and white regions have no interior equilibria in Figures  \ref{fig_InteriorModel1_rho2S_1} and \ref{fig_InteriorModel1_rho2S_2}.}
\label{fig3:M1Int3}
\end{figure}

\end{enumerate}

\begin{figure}[H]
 \vspace{-10pt}
\begin{center}
\subfigure[Time series of Model \ref{2DPatch} when $a_1 = 1$, $d_1 = 0.85$, $s = 0.55$, $\rho_1=13$, $x_1(0) = 1$, $y_1(0) = 0.25$, $ x_2(0) = 0.3$, and $y_2(0) = 0.7$.]{\includegraphics[height = 50mm, width =60mm]{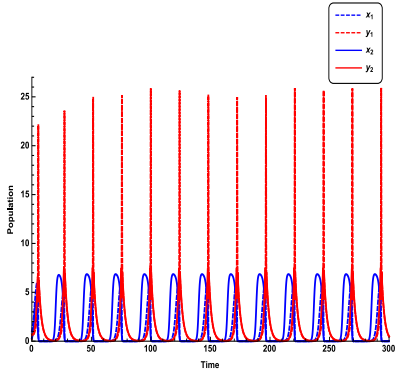}\label{fig_TimeSerie51}}\hspace{15mm}
\subfigure[Time series of Model \ref{2DPatch} when $a_1 = 2.1$, $d_1 = 2$, $\rho_1= 1$, $\rho_2=0.75$, $s = 0.85$, $x_1(0) = 0.9$, $y_1(0) = 1.1$, $ x_2(0) = 0.4$, and $y_2(0) = 0.8$]{\includegraphics[height = 50mm, width = 60mm]{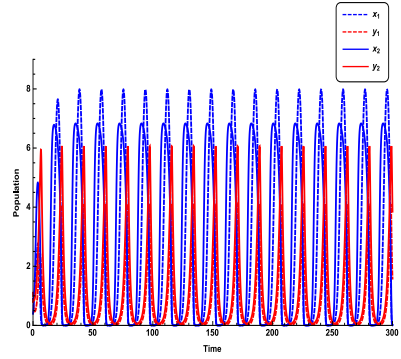}\label{fig_TimeSerie52}}
  \end{center}
\vspace{-15pt}
\caption{{\small Time series of Model \ref{2DPatch} when $r=1.8$, $d_2=0.35$, $K_1=10$, $K_2=7$, and $a_2=1.4$. Figures \ref{fig_TimeSerie51} and \ref{fig_TimeSerie52} illustrate the coexistence of prey and predator through fluctuating dynamics while Model \ref{2DPatch} has no interior equilibria. The blue dashed lines represent the prey population in patch 1, the dashed red lines represent the predator population in patch 1, the blue solid lines is the the prey population in patch 2, and the red solid lines represent predator population in patch 2.}} 
\label{fig:TimeSerie21}
\end{figure}

\noindent\textbf{No interior equilibrium but all species coexist with fluctuating dynamics:}  Our discussions above suggest that the large values of $s$ can destroy the interior equilibrium (see white regions in Figures \ref{fig_InteriorModel1_rho1S_1}, \ref{fig_InteriorModel1_rho1S_2}, \ref{fig_InteriorModel1_rho2S_1} and \ref{fig_InteriorModel1_rho2S_2}). Thus, the system is not permanent based on the fixed point theorem. However, our time series (e.g., Figures \ref{fig_TimeSerie51} and \ref{fig_TimeSerie52}) suggest that for almost all strictly positive initial conditions, both prey and predator can coexist through fluctuating dynamics for some white regions of Figures \ref{fig_InteriorModel1_rho2S_1} and \ref{fig_InteriorModel1_rho2S_2}.  \\

The proportion of the predators population engaging in the passive dispersal, i.e., $s$, has profound impacts on the population dynamics of prey and predator presented by Model \eqref{2DPatch} which generate complicated dynamics including different types of multiple attractors.\\

\begin{figure}[H]
 \vspace{-10pt}
\begin{center}
\subfigure[Time series of Model \ref{2DPatch} when $a_1 = 1$, $d_1 = 0.85$, $s = 0.8$, $x_1(0) = 0.05$, $y_1(0) = 1$, $ x_2(0) = 3.55$, and $y_2(0) = 2.7$ which converges to the boundary equilibrium $(x_1,y_1,x_2,y_2)=(0,1,3.6,2.9)$.]{\includegraphics[height = 57mm, width =67mm]{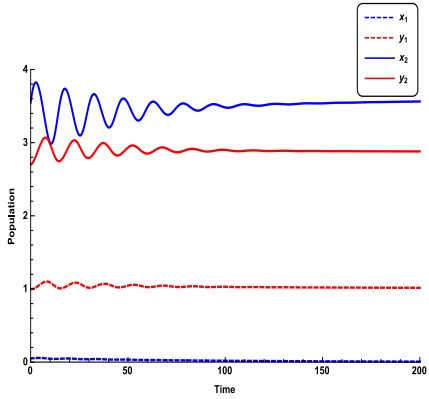}\label{fig_TimeSerie31}}\hspace{15mm}
\subfigure[Time series of Model \ref{2DPatch} when $a_1 = 1$, $d_1 = 0.85$, $s = 0.8$, $x_1(0) = 0.2$, $y_1(0) = 1.15$, $ x_2(0) = 2.7$, and $y_2(0) = 2.8$]{\includegraphics[height = 57mm, width = 67mm]{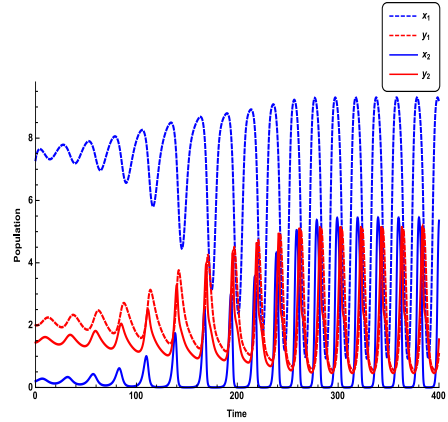}\label{fig_TimeSerie32}}
  \end{center}
\vspace{-15pt}
\caption{{\small Time series of Model \ref{2DPatch} when $r=1.8$, $d_2=0.35$, $K_1=10$, $K_2=7$, $a_2=1.4$, $\rho_1 = 1$, and $\rho_2 = 2.5$. Figures \ref{fig_TimeSerie31} and \ref{fig_TimeSerie32} represent the dynamical pattern generated by two interior saddles, one boundary sink and one boundary saddle. The blue dashed lines represent the prey population in patch 1, the dashed red lines represent the predator population in patch 1, the blue solid lines is the the prey population in patch 2, and the red solid lines represent predator population in patch 2.}} 
\label{fig:TimeSerie21}
\end{figure}

\noindent\textbf{Boundary attractor versus an interior attractor through two interior equilibria:} When Model \eqref{2DPatch} has two interior equilibria, the typical dynamics are that Model \eqref{2DPatch} either converges to a boundary attractor or  the interior attractor depending initial conditions. We provide an example in  Figures \ref{fig_TimeSerie31}, and \ref{fig_TimeSerie32} where $a_1 = 1$, $d_1 = 0.85$, $s = 0.8$ and $$r=1.8,\,\,d_2=0.35,\,\,K_1=10,\,\,K_2=7,\,\,a_2=1.4,\,\,\rho_1 = 1, \,\,\rho_2 = 2.5.$$ \\





\begin{figure}[H]
 \vspace{-10pt}
\begin{center} 
   \subfigure[Time series of Model \ref{2DPatch} when $a_1 = 1$, $d_1 = 0.85$, $s = 0.8392$, $x_1(0) = 0.25$, $y_1(0) = 1.05$, $ x_2(0) = 4.18$, and $y_2(0) = 2.68$ which stabilize at $(x_1,y_1,x_2,y_2)=(0.09,1.08,4.27,2.64).$]{\includegraphics[height = 57mm, width =67mm]{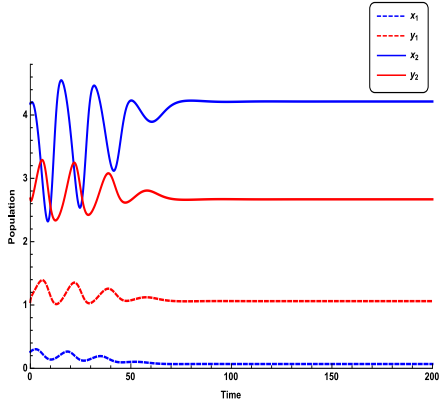}\label{fig_TimeSerie41}}\hspace{15mm}
   \subfigure[Time series of Model \ref{2DPatch} when $a_1 = 1$, $d_1 = 0.85$, $s = 0.8392$, $x_1(0) = 0.58$, $y_1(0) = 1.4$, $ x_2(0) = 2.5$, and $y_2(0) = 3.1$]{\includegraphics[height = 55mm, width =65mm]{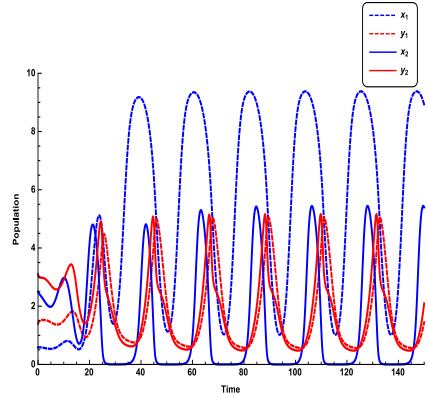}\label{fig_TimeSerie42}}
\end{center}
\vspace{-15pt}
\caption{{\small Time series of Model \ref{2DPatch} when $r=1.8$, $d_2=0.35$, $K_1=10$, $K_2=7$, $a_2=1.4$, $\rho_1 = 1$, and $\rho_2 = 2.5$. Figures \ref{fig_TimeSerie41} and \ref{fig_TimeSerie42} describe the dynamical pattern generated by two interior saddles and one interior that is locally stable. The blue dashed lines represent the prey population in patch 1, the dashed red lines represent the predator population in patch 1, the blue solid lines is the the prey population in patch 2, and the red solid lines represent predator population in patch 2.}} 
\label{fig:TimeSerie22}
\end{figure}

\noindent\textbf{Two interior attractors through three interior equilibria:} When Model \eqref{2DPatch} has three interior equilibria, the typical dynamics are that Model \eqref{2DPatch} has two interior attractors. We provide an example in Figures \ref{fig_TimeSerie41}, and \ref{fig_TimeSerie42} where $a_1 = 1$, $d_1 = 0.85$, $s = 0.8392$ and $$r=1.8,\,\,d_2=0.35,\,\,K_1=10,\,\,K_2=7,\,\,a_2=1.4,\,\,\rho_1 = 1, \,\,\rho_2 = 2.5.$$ \\

\section{Conclusion}\label{sec_discussion}

We propose and study a two patch prey predator model with the following assumptions: (1) Only predators can migrate and preys are immobile; (2) predators use two dispersal strategies: the passive dispersal and the predation attraction; (3) The model is reduced to the Rosenzweig-MacArthur model in the absence of dispersal. We provide boundedness and positivity of the proposed model in Theorem \ref{th1:pb}. The analytical results which is summarize in Table \eqref{table_stability} along with the numerical results presented throughout the paper answer the questions regarding the dynamics of our proposed nonlinear model:\\

 When there is no prey in one of the patches, our model applies to the sink-source dynamics where no prey patch is the sink. Analytical results (Theorem \ref{th4:persist}) imply that  predators could be driven to extinction locally if the product of the dispersal strength and the proportion of predator population using the passive dispersal (i.e. $s$) are large. In addition, the sink-source dynamics can process two interior equilibria (see Proposition \eqref{pr1:be}). Our simulations (Figure \ref{fig:intm3_1}) suggest that the small values of $s$ lead to permanence of the system which is supported by Theorem \ref{th4:persist}. For the intermediate values of $s$, the system can can have two interiors $E_{x_1,y_1,y_2}^{l},l=1,2$ ($i=1$) or $E_{y_1,x_2,y_2}^{l},l=1,2$ ($i=2$); For the large values of $s$, it has no interior equilibria with the consequences that predator goes extinct in two patches. In addition, the intermediate values  of $s$ can stabilize the dynamics with certain dispersal strengths (see blue line for locally stable in Figures \ref{fig_InteriorModel3x20_1}, \ref{fig_InteriorModel3x10_1}, and \ref{fig_InteriorModel3x10_2}).  \\
 
 Theorem \eqref{th2:be} and Proposition \eqref{pr1:be} provide the existence of the boundary equilibria and the related local stability of  our model \eqref{2DPatch}. These results illustrate how $s$ can potentially stabilize the basic boundary equilibria $E_{K_10K_20}$ consequently driving predator extinct in both patches locally.  Theorem \eqref{th34:interior} provide insights into the existence and stability of a symmetric interior equilibria when Model \eqref{2DPatch} is symmetric (i.e. in exception of the dispersal strength and dispersal strategy, all life history parameters are the same in both patches). The analytical results indicate  that the dispersal strategies do not affect the existence and stability of this symmetric interior equilibria denoted $E$. However, bifurcation diagrams shown in Figures \ref{fig_InteriorSyModel1_Sy1} and \ref{fig_InteriorSyModel1_rho1S} suggest that the large predator population using the passive dispersal could generate two additional asymmetric interior equilibria which can be saddle or locally stable, thus generate bistability between two different attractors (see blue lines in Figure \ref{fig_InteriorSyModel1_Sy1} when $0.78\leq s \leq 0.92$). \\

Our numerical simulations performed in Section 4 show that the dispersal strategies, i.e., the portion of predator population using the passive dispersal strategies, have huge impacts on the prey and predator populations in two patches. The intermediate predator population using the passive dispersal  tends to stabilize the dynamics. Depending on the other life history parameters, the large or small predator population using the passive dispersal with certain dispersal strengths could generate multiple interior equilibria (up to three interior equilibria), thus lead to multiple attractors potentially. When Model \eqref{2DPatch} has two interior equilibria, it either converges to a boundary attractor or  the interior attractor depending initial conditions (see Figures \ref{fig_TimeSerie31}, and \ref{fig_TimeSerie32}); when Model \eqref{2DPatch} has three interior equilibria, it can have two interior equilibria (see Figures \ref{fig_TimeSerie41}, and \ref{fig_TimeSerie42}). The large predator population using the passive dispersal combined with the large dispersal strength can destroy the interior equilibria with consequences that prey in one patch may go extinct but predator persists in each patch. However, there are situations when the two patch model has no interior equilibrium but all species coexist with fluctuating dynamics (see Figures \ref{fig_TimeSerie51} and \ref{fig_TimeSerie52}). \\

The summary of our finding illustrates how population dynamics of prey and predators are affected by changing their foraging behavior. This study give us a better understanding on how combinations of different foraging strategies used by predator favor or affect their coexistence or extinction. Many species tend to adapt to environmental conditions and change their foraging behavior accordingly (see example of foraging behavior of Ants in \citet{taylor1977foraging, markin1970foraging, traniello1984ant}). It will be interesting to look at a two patch prey predator model with adaptive foraging behavior in which adaptation is driven by certain environmental conditions such as temperature or availability of local resources. Such work is on going by the authors.


\section{Appendix: Proofs}

 \subsection*{ Proof of Theorem \ref{th1:pb}}
 \begin{proof}
Observed that $\frac{dx_i}{dt} \big\vert_{x_i=0}=0$ and $\frac{dy_i}{dt} \big\vert_{y_i=0}=\rho_i s y_j \geq 0$ if $ y_j \geq 0$ for $i = 1, 2$, $j = 1,2$, and $i \not=j$. The model \eqref{2DPatch} is positively invariant in $\mathbb R^4_+$ by theorem A.4 (p. 423) in \citet{thieme2003mathematics}. It follows that the set $\{(x_1,y_1,x_2,y_2)\in\mathbb R^4_+:x_i=0 \}$ is invariant for both $i = 1,2$ under the same theorem. The proof of boundedness is as follow
$$\frac{dx_i}{dt}= r_i x_i\left(1 - \frac{x_i}{K_i}\right) - \frac{a_i x_i y_i}{1 + x_i}\leq r_i x_i\left(1 - \frac{x_i}{K_i}\right)$$ thus $$\limsup_{t\rightarrow\infty} x_i(t)\leq K_i.$$
Now we define $L = \rho_2 (x_1+y_1)+\rho_1 (x_2+y_2)$, to get
\begin{align*}
\frac{dL}{dt}&=\rho_2\frac{d(x_1+y_1)}{dt}+\rho_1\frac{d(x_2+y_2)}{dt}\\
&=\rho_2x_1\left(1 - \frac{x_1}{K_1}\right) +\rho_1 rx_2\left(1 - \frac{x_2}{K_2}\right)  -\rho_2 d_1 y_1 -\rho_1 d_2 y_2 + \rho_1\rho_2(1-s)\left( \frac{a_1 x_1 y_1y_2}{1 + x_1}-\frac{a_2 x_2 y_2y_1}{1 + x_2}\right) \\
&- \rho_1\rho_2(1-s)\left( \frac{a_1 x_1 y_1y_2}{1 + x_1}-\frac{a_2 x_2 y_2y_1}{1 + x_2}\right)+\rho_1\rho_2s(y_1-y_2) - \rho_1\rho_2s(y_1-y_2)\\
&=\rho_2x_1\left(1 - \frac{x_1}{K_1}+d_1\right) +\rho_1 x_2\left(r - \frac{rx_2}{K_2}+d_2\right) -\rho_2 d_1 (x_1+y_1) -\rho_1 d_2 (x_2+y_2)\\
&\leq T- d_{\text{min}} \left[\rho_2 (x_1+y_1)+\rho_1 (x_2+y_2)\right]=T- d_{\text{min}}L
\end{align*}
 where $d_{\text{min}}=\min\{d_1,d_2\}$ and 
$$T=\max_{0\leq x_1\leq K_1}\Big\{\rho_2x_1\left(1 - \frac{x_1}{K_1}+ d_1\right)\Big\} +\max_{0\leq x_2\leq K_2}\Big\{\rho_1x_2\left(r - \frac{rx_2}{K_2}+d_2\right)\Big\}.$$
consequently
$$\limsup_{t\rightarrow\infty} L(t)=\limsup_{t\rightarrow\infty} \rho_2 (x_1(t)+y_1(t))+\rho_1 (x_2(t)+y_2(t))\leq \frac{T}{d_{\text{min}}}.$$ This shows that Model \eqref{2DPatch} is bounded in $\mathbb R^4_+$ which conclude the proof of theorem \eqref{th1:pb}.
\end{proof}

 \subsection*{ Proof of Theorem \ref{th4:persist}}
 \begin{proof}

Item 1: Model \eqref{2DPatch} is positively invariant and bounded in $\mathbb R^4_+$ according to Theorem \ref{th1:pb}. From this, it follows that Model \eqref{2DPatch} is attracted to a compact set $C$ in $\mathbb R^4_+$. Furthermore, if $x_j=0,~j=1,\mbox{  or  }2$, then Model \eqref{2DPatch} is reduced to three species couple models \eqref{xi0}. Consider the fact that $\lim_{t\rightarrow\infty} y_i(t)=\lim_{t\rightarrow\infty} y_j(t)=0$ when $x_i = 0$, we can conlcude that $y_1 =y_2 = 0$ is an omega limit set of Model \eqref{xi0}. Additionally 
$$\frac{dx_i}{x_idt}\Big\vert_{x_i=0}=r_i>0$$ 
then by Theorem 2.5 of Hutson (1984) \citep{hutson1984theorem}, prey $x_i$ persists.\\

Item 2: Define $V(y_i,y_j)=\rho_jy_i+\rho_iy_j$, then we have
\begin{align*}
\frac{dV}{dt}=\frac{ax_iy_i}{1+x_i}\rho_j-d_iy_i\rho_j-d_jy_j\rho_i =  \left[\frac{a_ix_i}{1+x_i}-d_i\right]y_i\rho_j-d_jy_j\rho_i.
\end{align*}
Notice that $\limsup_{t\rightarrow\infty} x_i(t)\leq K_i$. Then if $\mu_i>K_i$ we have $\frac{a_iK_i}{1+K_i}-d_i<-\delta<0$ and let $\delta^*=\min\{\delta,d_j\}$. This implies 
\begin{align*}
\frac{dV}{dt} =  \left[\frac{a_ix_i}{1+x_i}-d_i\right]y_i\rho_j-d_jy_j\rho_i &<-(\delta y_i\rho_j+d_jy_j\rho_i) \\
&< -\delta^* ( y_i\rho_j+d_jy_j\rho_i)=-\delta^* V(y_i,y_j)
\end{align*}
Therefore both predators go extinct if $\mu_i>K_i$. Now Model \eqref{xi0} reduces to the following prey model since $\limsup_{t\rightarrow\infty} y_i(t) = \limsup_{t\rightarrow\infty} y_j(t) =0$
$$\frac{dx_i}{dt} = r_ix_i\left(1-\frac{x_i}{K_i}\right)$$
with prey $x_i$ converging to $K_i$. Thus Model \eqref{xi0}is globally stable at $(K_i,0,0)$ when $\mu_i>K_i$.\\

Item 3: Now we focus on the persistence condition for predator $y_i$. Since $x_i$ is persistent from Item 1 Theorem  \ref{th4:persist} then we can conclude that Model \eqref{xi0} is attracted to a compact set $C_s$ subset of $C$ that exclude $E_{000}$.Then according to Theorem \ref{th1:pb} and \ref{th2:be}, the omega limit set of Model \eqref{xi0} on the compact set $C_s$ is $E_{K_i00}$. 
Notice that the following inequalities,
\begin{align*}
\frac{dy_i}{dt} & =  \frac{a_i x_i y_i}{1 + x_i} -d_i y_i+\rho_i(1-s)\left(\frac{a_ix_iy_i}{1+x_i}y_j\right)+ \rho_is( y_j-y_i)\\
&\geq \frac{a_ix_iy_i}{1+x_i}-d_iy_i-\rho_isy_i \quad \Rightarrow \\
\quad \frac{dy_i}{y_idt} &\geq \frac{a_ix_i}{1+x_i}-d_i-\rho_is\\
\end{align*}
therefore, we have 
$$\frac{dy_i}{y_idt} \Big\vert_{E_{K_i00}}  \geq \frac{a_iK_i}{1+K_i} - (d_i+\rho_is).$$
According to Theorem 2.5 of Hutson (1984) \citep{hutson1984theorem}, we can conclude that predator $y_i$ is persistent if the following inequalities hold
$$\frac{a_iK_i}{1+K_i} - (d_i+\rho_is)>0\Leftrightarrow \rho_is<\frac{(a_i-d_i)(K_i-\mu_i)}{1+K_i}.$$
Now assume that $\rho_is<\frac{(a_i-d_i)(K_i-\mu_i)}{1+K_i}$ holds, then we can conclude that predator $y_i$ is persistent. This implies that when time large enough, there exists some $\epsilon>0$ such that 
$$\frac{dy_i}{dt}\big\vert_{y_i=0}=\rho_j s y_j>\rho_j s \epsilon>0.$$
Thus, we could conclude that predator in Patch $j$ also persists due to the persistence of predator in Patch $i$.\\

 \end{proof}
 

 \subsection*{ Proof of Proposition \ref{pr1:be}}
 \begin{proof}
 The algebraic calculations imply that an interior equilibrium $(x_i^*,y_i^*,y_j^*)$ of Model \eqref{xi0} satisfies the following equations:
 $$\begin{aligned}
 y_{i}^* &=\frac{r_i(K_i-x_i^*)(1+x_i^*)}{a_iK_i} \\
 y_j^* &= \frac{r_i(K_i-x_i^*)[x_i^*(a_i-d_i)-d_i]\rho_j}{a_iK_id_j\rho_i}\\
 0&=\underbrace{[(x_i^*)^{3}-(\mu_i+K_i)(x_i^*)^{2}-\alpha_i x_i^*+\beta_i]}_{f_i(x_i^*)}[x_i^*+1] \\
 \end{aligned}$$ where $\beta_i = \frac{\left[d_j\rho_is+d_i(d_j+\rho_js)\right]K_i}{r_i(a_i-d_i)(1-s)\rho_j}$ and
  {\footnotesize $$\alpha_i =\frac{\left[d_js\rho_i+r_id_i(1-s)-(a_i-d_i)(d_j+s\rho_j)\right]K_i}{r_i(a_i-d_i)(1-s)\rho_j} =\beta_i +\frac{[r_id_i(1-s)-a_i(d_j+s\rho_j)]K_i}{r_i(a_i-d_i)(1-s)\rho_j} .$$}
 This implies that
 $$ 0<\mu_i=\frac{d_i}{a_i-d_i}<x_i^*<K_i \mbox{ and } f_i(x_i^*)=0.$$
 Therefore, if $a_i<d_i$ or $\mu_i>K_i$ or $f_i(x_i^*)=0$ has no positive roots, then Model \eqref{xi0} has no interior equilibrium.\\
 
 Now assume that $ 0<\mu_i=\frac{d_i}{a_i-d_i}<K_i$, then we have $f_i(0) = \beta_i >0$ and $\lim_{x_i\rightarrow -\infty}f_i(x_i) =-\infty$. This indicates that there exist $x_0 \in(-\infty,0)$ such that $f_i(x_0)= 0$. Therefore, we can conclude that $f_i(x_i)$ has at least one negative root and at most two positive roots since  $f_i(x_i)$ is a polynomial with degree 3. The derivative of $f_i(x_i)$ has the following form
 {\footnotesize
 \begin{align*}
 f_i^{'}(x_i) = 3x_i^{2}-2(\mu_i+K_i)x_i-\alpha_i = 0 
 \end{align*}
 } which gives the following two critical points if $\Delta_i = (\mu_i+K_i)^2+3\alpha_i>0$
$$  x_{i}^{c_{+,-}} = \frac{(\mu_i+K_i)\pm \sqrt{(\mu_i+K_i)^2+3\alpha_i}}{3} = \frac{(\mu_i+K_i)\pm\sqrt{\Delta}}{3}.$$
Therefore if $\Delta_i \geq 0$, then $x_i^{c_{+}}= \frac{(\mu_i+K_i)+ \sqrt{\Delta_i}}{3} > 0$ is the local minimum of $f_i(x_i)$ since $f_i^{''}(x_i^{c_{+}}) = 2 \sqrt{\Delta_i}\geq 0$ and $f_i^{''}(x_i^{c_{-}}) = -2 \sqrt{\Delta_i}\leq 0$. We note that $f_i(x_i)$ has two positive roots if $f_i(x_i^{c_{+}}) \leq 0$. It follows that $f_i(x_i)$ has two positive roots if the following equation is satisfied:

 {\footnotesize
 \begin{align*}
 f\left(x_i^{c_{+}}\right) =-\frac{1}{3}\left[\alpha_i(\mu_i+K_i)-3\beta_i\right]  -\frac{1}{27}\left[(\mu_i+K_i)+3\alpha_i\right]^2\left[2(\mu_i+K_i)-\sqrt{\Delta_i}\right]-\frac{1}{3}\alpha_i\sqrt{\Delta_i} < 0.
 \end{align*}
 }
Since 
 $$\alpha_i(\mu_i+K_i)-3\beta_i > 0 \quad \Rightarrow \quad \alpha_i>\frac{3\beta_i}{\mu_i+K_i}$$
 and
$$2(\mu_i+K_i)-\sqrt{\Delta_i} = 2(\mu_i+K_i)-\sqrt{(\mu_i+K_i)^2+3\alpha_i} > 0 \quad \Rightarrow \quad \alpha_i < (\mu_i+K_i)^2$$
  therefore we can conclude that $f_i(x_i)$ has two positive roots when $\frac{3\beta_i}{\mu_i+K_i}<\alpha_i<(\mu_i+K_i)^2$.  Thus for $x_{i\ell}^*$ where $\ell=1,2$ denote the two positive roots of the nullclines $f_i(x_i)$ and $i=1,2$ represent the prey population in patch one and two, we have:
{\footnotesize
\begin{align*}
 y_{i\ell}^* =\frac{(K_i-x_{i\ell}^*)(1+x_{i\ell}^*)}{a_iK_i}, \hspace{1in}  y_{j\ell}^* =\frac{(K_i-x_{i\ell}^*)[x_{i\ell}^*(a_i-d_i)-d_i]\rho_j}{a_iK_id_j\rho_i}
 \end{align*} 
 } if $\mu_i<x_{i\ell}^*<K_i, \ell=1,2.$
 
 From the arguments above we conclude that Model \eqref{xi0} can have up to two interior equilibria $E_{x_i,y_i,y_j}^{\ell}=(x_{i\ell}^*,y_{i\ell}^*,y_{j\ell}^*)$ when  $\frac{3\beta_i}{\mu_i+K_i}<\alpha_i<(\mu_i+K_i)^2$ and $\mu_i<x_{i\ell}^*<K_i, \ell=1,2.$\\

 On the other hand, if $\Delta_i=  (\mu_i+K_i)^2+3\alpha_i< 0$ then $f_i(x_i) $ has no positive real roots and hence Model \eqref{xi0} has no interior equilibrium.\\
\end{proof}


 
 \subsection*{ Proof of Theorem \ref{th2:be}}
\begin{proof}

The local stability of the equilibrium $(x_1^*,y_1^*,x_2^*,y_2^*)$  of Model \eqref{2DPatch} is established by finding the eigenvalues $\lambda_i$, $i = 1,~2,~3,~4$ of the Jacobian matrix $J_{(x_1^*,y_1^*,x_2^*,y_2^*)}$ \eqref{JE} evaluated at the equilibria. 
 
 \begin{center}
{\miniscule
\bae\label{JE}
\begin{array}{ll}
&J_{(x_1^*,y_1^*,x_2^*,y_2^*)}=\\\\
&\left[\begin{array}{cccc}
\lx1-\frac{2x_1^*}{K_1}\rx-\frac{a_1y_1^*}{\lx1+x_1^*\rx^2} & -\frac{a_1x_1^*}{1+x_1^*} & 0 & 0\\

\frac{a_1y_1^*\lx 1+y_2^*\lx \rho_1-s\rho_1\rx\rx}{\lx1+x_1^*\rx^2} & \rho_1\lx1-s\rx\lx\frac{a_1x_1^*y_2^*}{1+x_1^*}-\frac{a_2x_2^*y_2^*}{1+x_2^*}\rx + \frac{a_1x_1^*}{1+x_1^*}-d_1-s\rho_1 & \frac{\rho_1 a_2\lx-1+s\rx y_1^*y_2^*}{(1+x_2^*)^2} & s\rho_1+\rho_1\lx 1-s\rx\lx\frac{a_1x_1^*y_1^*}{1+x_1^*} - \frac{a_2x_2^*y_1^*}{1+x_2^*}\rx \\

0 & 0 & r\lx1-\frac{2x_2^*}{K_2}\rx-\frac{a_2y_2^*}{\lx1+x_2^*\rx^2} & -\frac{a_2x_2^*}{1+x_2^*} \\

\frac{\rho_2 a_1\lx-1+s\rx y_1^*y_2^*}{(1+x_2^*)^2} & s\rho_2+\rho_2\lx 1-s\rx\lx\frac{a_2x_2^*y_2^*}{1+x_2^*} - \frac{a_1x_1^*y_1^*}{1+x_1^*}\rx & \frac{a_2y_2^*\lx 1+y_1^*\lx \rho_2-s\rho_2\rx\rx}{\lx1+x_2^*\rx^2} &  \rho_2\lx1-s\rx\lx\frac{a_2x_2^*y_1^*}{1+x_2^*}-\frac{a_1x_1^*y_1^*}{1+x_1^*}\rx + \frac{a_2x_2^*}{1+x_2^*} - d_2-s\rho_2 \\
\end{array}\right]
\end{array}
\eae}
 \end{center}
 
By substituting the equilibria $E_{0000}, E_{K_1000}, E_{00K_20,}$ into the Jacobian matrix \eqref{JE}, it was found that these equilibria are saddle consider one of their eigenvalues is positive. \\
For the equilibrium $E_{K_10K_20}$ we obtain
$$\lambda_1 = -1~(<0),~~\lambda_2 = -r~(<0),$$

 $$\lambda_3+\lambda_4= \frac{a_1K_1}{1+K_1}-d_1 + \frac{a_2K_2}{1+K_2}-d_2 - s\rho_1 - s\rho_2  $$
 and
 $$\lambda_3\lambda_4 = \left[d_1-\frac{a_1K_1}{1+K_1}\right]\left[1-\frac{a_2K_2}{(s\rho_2+d_2)(1+K_2)}\right] +  \frac{s\rho_1}{s\rho_2+d_2}\left[d_2-\frac{a_2K_2}{1+K_2}\right]$$
 
 Notice that the eigenvalue $\lambda_3$ and $\lambda_4$ being negative for $s \in(0, 1)$ is equivalent to the case where the boundary equilibria $(K_i,0)$ for the single patch is globally asymptotically stable. This is also equivalent to $\mu_i > K_i$ or $\frac{a_iK_i}{1+K_i} - d_i < 0$. We again observe that for $\frac{a_iK_i}{1+K_i} - d_i < 0$ the following holds
 $$\lambda_3+\lambda_4 =  \frac{a_1K_1}{1+K_1}-d_1 + \frac{a_2K_2}{1+K_2}-d_2 - s\rho_1 - s\rho_2 < 0\Rightarrow d_1+ d_2+ s\rho_1 + s\rho_2 > \frac{a_1K_1}{1+K_1} + \frac{a_2K_2}{1+K_2}$$
and
$$\lambda_3\lambda_4 = \left[d_1-\frac{a_1K_1}{1+K_1}\right]\left[1-\frac{a_2K_2}{(s\rho_2+d_2)(1+K_2)}\right] +  \frac{s\rho_1}{s\rho_2+d_2}\left[d_2-\frac{a_2K_2}{1+K_2}\right] > 0 $$
which can be rewritten in the following form:

\begin{align*}
& \mathlarger{\sum}_{i=1}^{2}\left[ \frac{(a_i-d_i)(\mu_i-K_i)}{1+K_i} + s\rho_i \right] > 0 \\
&\mbox{  and  } \\
&\left[\frac{(a_1-d_1)(\mu_1-K_1)}{1+K_1}\right]\left[s\rho_2+\frac{(a_2-d_2)(\mu_2-K_2)}{1+K_2}\right]+s\rho_1\left[\frac{(a_2-d_2)(\mu_2-K_2)}{1+K_2}\right] > 0.
\end{align*} 

 Based on the discussion above, we can conclude that the results on the local stability of four boundary equilibria of Theorem \ref{th2:be} holds.\\
\\


Item 1:  Let $p_i(x)=\frac{a_i x}{1+x}$ and $q_i(x)=\frac{r_i(K_i-x)(1+x)}{a_iK_i}$ then we have the following
\begin{align*}
\frac{dx_i}{dt} &=r_ix_i\left(1-\frac{x_i}{K_i}\right)-\frac{a_i x_iy_i}{(1+x_i)}=\frac{a_i x_i}{1+x_i}\left[\frac{r_i(K_i-x_i)(1+x_i)}{a_iK_i}-y_i\right]=p_i(x_i)\left[q_i(x_i)-y_i\right]. \\
\frac{dy_i}{dt} &= y_i\left[\frac{a_ix_i}{1+x_i} - d_i\right] + \rho_i(1-s)y_iy_j\left[\frac{a_ix_i}{1+x_i} - \frac{a_jx_j}{1+x_j}\right]+\rho_is\left[y_j-y_i\right] \\
& = y_i\left[p_i(x_i) - d_i\right] + \rho_i(1-s)y_iy_j\left[p_i(x_i) - p_j(x_j)\right]+\rho_is\left[y_j-y_i\right] \mbox{ where both } i, j = 1,~2, \mbox{ with } i\not=j
\end{align*}
Now consider the following Lyapunov functions
\begin{equation}\label{Lyapunov1}
\begin{aligned}
V_1(x_1,y_1) =\rho_2 \int^{x_1}_{K_1} \frac{p_1(\xi)-p_1(K_1)}{p_1(\xi)}d\xi + \rho_2y_1
\end{aligned}
\end{equation}
and
\begin{equation}\label{Lyapunov2}
\begin{aligned}
V_2(x_2,y_2) =\rho_1 \int^{x_2}_{K_2} \frac{p_2(\xi)-p_2(K_2)}{p_2(\xi)}d\xi + \rho_1y_2
\end{aligned}
\end{equation}
Taking derivative of the functions \eqref{Lyapunov1} and \eqref{Lyapunov2} with respect to time $t$ yield
{\footnotesize
\begin{equation}\label{L1}
\begin{aligned}
\frac{d}{dt}V_1(x_1(t),y_1(t)) &= \rho_2\frac{p_1(x_1)-p_1(K_1)}{p_i(x_1)}\frac{dx_1}{dt} + \rho_2\frac{dy_1}{dt} \\
 &= \rho_2\frac{1}{p_1(x_1)}\lz p_1(x_1)-p_1(K_1)\rz p_1(x_1)\lz q_1(x_1)-y_1\rz +\rho_2 y_1\lz p_1(x_1)-d_1\rz + \rho_1\rho_2(1-s) y_1y_2\lz p_1(x_1)-p_2(x_2)\rz \\
 &+ \rho_1\rho_2 s\lz y_2 - y_1\rz\\
 &=\rho_2\lz p_1(x_1)-p_1(K_1)\rz q_1(x_1)+ \rho_2 y_1\lz p_1(K_1)-d_1\rz+\rho_1\rho_2(1-s) y_1y_2\lz p_1(x_1)-p_2(x_2)\rz + \rho_1\rho_2 s\lz y_2 - y_1\rz
\end{aligned}
\end{equation}
}
and
{\footnotesize
\begin{equation}\label{L2}
\begin{aligned}
\frac{d}{dt}V_2(x_2(t),y_2(t)) &=\rho_1 \frac{p_2(x_2)-p_2(K_2)}{p_2(x_2)}\frac{dx_2}{dt} +\rho_1 \frac{dy_2}{dt} \\
&= \rho_1\frac{1}{p_2(x_2)}\lz p_2(x_2)-p_2(K_2)\rz p_2(x_2)\lz q_2(x_2)-y_2\rz +\rho_1 y_2\lz p_2(x_2)-d_2\rz + \rho_1\rho_2(1-s) y_1y_2\lz p_2(x_2)-p_1(x_1)\rz \\
&+ \rho_1\rho_2s\lz y_1 -y_2\rz\\
&=\rho_1\lz p_2(x_2)-p_2(K_2)\rz q_2(x_2)+ \rho_1 y_2\lz p_2(K_2)-d_2\rz - \rho_1\rho_2(1-s) y_1y_2\lz p_1(x_1)-p_2(x_2)\rz - \rho_1\rho_2s\lz y_2 -y_1\rz
\end{aligned}
\end{equation}
}
Also, we denote $V=V_1+V_2$ and adding \eqref{L1} and \eqref{L2}, we obtain
{\footnotesize
\begin{equation}\nonumber
\begin{aligned}
\frac{d}{dt}V &=\frac{d}{dt}V_1(x_1(t),y_1(t)) + \frac{d}{dt}V_2(x_2(t),y_2(t)) \\
& = \rho_2\lz p_1(x_1)-p_1(K_1)\rz\lz q_1(x_1)-y_1\rz + \rho_2y_1\lz p_1(x_1)-d_1\rz + \rho_1\lz p_2(x_2)-p_2(K_2)\rz\lz q_2(x_2)-y_2\rz + \rho_1y_2\lz p_2(x_2)-d_2\rz \\
&=  \rho_2\lz p_1(x_1)-p_1(K_1)\rz q_1(x_1)+ \rho_2 y_1\lz p_1(K_1)-d_1\rz + \rho_1\lz p_2(x_2)-p_2(K_2)\rz q_2(x_2)+ \rho_1 y_2\lz p_2(K_2)-d_2\rz.
\end{aligned}
\end{equation}
}
We observe that the function $p_i(x_i)$ increases as $x_i$ increases thus $p_i(x_i) - p_i(K_i) > 0$ if $x_i > K_i$ and $p_i(x_i) - p_i(K_i) < 0$ if $x_i < K_i$. Also, $q_i(x_i)$ is positive if $x_i < K_i$ and it is negative if $x_i > K_i$. This implies that the expressions $ \rho_2\lz p_1(x_1)-p_1(K_1)\rz q_1(x_1)$ and $\rho_1\lz p_2(x_2)-p_2(K_2)\rz q_2(x_2)$ are both negative for all $x_i \geq 0$ since all the parameters are assumed to be positive. Also, Assume $\mu_i> K_i$. This implies that $\frac{d_i}{a_i-d_i}>K_i$ which is also equivalent to $\frac{a_iK_i}{1+K_i}=p_i(K_i)<d_i$. Since $p_i(K_i) < d_i$ then $p_i(K_i) - d_i < 0$. The derivative $\frac{dV}{dt}$ is therefore negative which implies that both $V_1$ and $V_2$ are Lyapunov functions, and the boundary equilibrium $E_{K_10K_20} = \lx K_1,0,K_2,0\rx$ is globally stable when $\mu_i> K_i$ by Theorem $3.2$ in \cite{hsu1978global}.\\

Item 2: According to Theorem \ref{th1:pb}, we know that Model \eqref{2DPatch} is attracted to a compact set $C$ in $\mathbb R^4_+$. Define $V_x=x_1+x_2$, then we have
$$\frac{dV_x}{dt}=\frac{dx_1}{dt}+\frac{dx_2}{dt}= r_1x_1\left(1 - \frac{x_1}{K_1}\right) - \frac{a_1 x_1 y_1}{1 + x_1}+r_2 x_2\left(1 - \frac{x_2}{K_2}\right) - \frac{a_2 x_2 y_2}{1 + x_2}.$$
Notice that if $x_i=x_j=0$, then Model \eqref{2DPatch} converges to $(0,0,0,0)$, and
$$\frac{dV_x}{dt}\Big\vert_{x_1=x_2=0}=r_1+r_2>0.$$
Therefore, according to Theorem 2.5 of \cite{hutson1984theorem}, we can conclude that prey population in two patches, i.e., $x_1+x_2$, is persistent. Moreover, if $x_j=0$, Model \eqref{2DPatch} is reduced to the subsystem \eqref{xi0} where prey $x_i$ is persistent according to Theorem \ref{th4:persist}. Thus, we can conclude prey population in at least one patch is persistent.\\

Define $V_y=\rho_2y_1+\rho_1y_2$, then we have
$$\frac{dV_y}{dt}=\rho_2\frac{dy_1}{dt}+\rho_1\frac{dy_2}{dt}=\rho_2 y_1\left( \frac{a_1 x_1 }{1 + x_1} -d_1 \right)+\rho_1 y_2\left( \frac{a_2 x_2 }{1 + x_2} -d_2 \right).$$
Notice that if $y_i=y_j=0$, then Model \eqref{2DPatch} converges to  $(K_1,0,K_2,0)$. Since we have $K_i>\mu_i$ for both $i=1,2$, then we have 
$$\min_{i=1,2}\{ \frac{a_i K_i }{1 + K_i} -d_i \}=\delta>0.$$This implies that
$$\frac{dV_y}{dt}\Big\vert_{y_1=y_2=0}=\rho_2 y_1\left( \frac{a_1 K_1 }{1 + K_1} -d_1 \right)+\rho_1 y_2\left( \frac{a_2 K_2 }{1 + K_2} -d_2 \right)\geq \delta (\rho_2 y_1+\rho_1 y_2)=\delta V_y>0.$$
Therefore, according to Theorem 2.5 of \cite{hutson1984theorem} and the proof of Proposition \ref{pr1:be}, we can conclude that predator population in each patch is persistent.\\
 
\end{proof}
 
 
  \subsection*{Proof of Theorem \ref{th34:interior}}
 \begin{proof}
First we show the existence of the interior equilibrium $E = (\mu,\nu,\mu,\nu)$ in the symmetric case (i.e. $a_1=a_2=a, d_1=d_2=d, K_1=K_2=K,r=1$). The interior equilibrium can be obtained by the positive intersection of the two nullclines $x_1=F_1(x_2)$ and $x_2=F_2(x_1)$  \eqref{interior-eq3}. Recall from the nullcines \eqref{interior-eq3} that
$$x_i(x_j) = \frac{(\mu+K)\pm\sqrt{(\mu+K)^2-4\phi_i(x_j)}}{2} .$$
where $ \phi_i(x_j)=\mu K+\frac{\rho_i}{\rho_j}(x_j-\mu)(x_j-K)$ which indicate that 

{\footnotesize
$$x_i^{+}(\mu) = \frac{(\mu+K)+\sqrt{(\mu+K)^2-4\phi_i(\mu)}}{2} = K \mbox{  and  } x_i^{-}(\mu) = \frac{(\mu+K)-\sqrt{(\mu+K)^2-4\phi_i(\mu)}}{2} = \mu$$
}
This implies that $x = \mu$ is a positive solution of the nullcline \eqref{interior-eq3} when $a > d$ in the symmetric case. We can accordingly say that $E=(\mu,\nu,\mu,\nu)$ is an interior equilibrium of Model \eqref{2DPatch} when $a_1=a_2=a, d_1=d_2=d, K_1=K_2=K,r=1$. \\

The local stability of $E=(\mu,\nu,\mu,\nu)$ is obtained by the eigenvalues of the Jacobian matrix \eqref{JE} evaluated at this equilibrium as follow:
{\footnotesize
\begin{align*}
\lambda_1\lambda_2 &= \frac{d(K-\mu)}{K(1+\mu)} > 0 \mbox{  if  } K > \mu \mbox{  and  } \lambda_1\lambda_2 = \frac{d(K-\mu)}{K(1+\mu)} < 0  \mbox{  if  } K < \mu\\[2mm]
\lambda_1 + \lambda_2 &= \frac{K-1-2\mu}{K(1+\mu)} < 0  \mbox{  if  } \mu > \frac{K-1}{2} \mbox{  and  } \lambda_1 + \lambda_2 = \frac{d(K-\mu)}{K(1+\mu)} > 0  \mbox{  if  }  \mu < \frac{K-1}{2}\\[2mm]
\lambda_3\lambda_4 &= \frac{(\rho_1+\rho_2)[(1-s)(K-\mu)d\nu-((K-1)-2\mu)s\mu)]+d(K-\mu)}{K(1+\mu)}> 0 \mbox{  for  } K > \mu \mbox{  and  } \mu > \frac{K-1}{2} \mbox{  when  } s \in [0,1] \\[2mm]
\lambda_3 + \lambda_4 &= -\left[\frac{-\mu(K-1)+2\mu^2+Ks(\rho_1+\rho_2)(1+\mu)}{K(1+\mu)}\right] < 0 \mbox{  for  } \mu > \frac{K-1}{2} \mbox{  when  } s \in [0,1]
\end{align*} 
}
Notice that the eigenvalues $\lambda_1,~\lambda_2,~\lambda_3,\mbox{  and  }\lambda_4$ being negative correspond to the case where the unique interior equilibrium $(\mu,\nu)$ of the single patch Model \eqref{Onepatch} is locally asymptoticaly stable. We can hence conclude that $E$ has the same local stability as the interior equilibrium $(\mu,\nu)$ for the single patch model \eqref{Onepatch}. Consequently $\frac{K-1}{2} < \mu <K$ are sufficients conditions for $E=(\mu,\nu,\mu,\nu)$ to be locally asymptotically stable while unstable when $\mu < \frac{K-1}{2} $ for $s \in [0,1]$.\\

\end{proof}

 
{
\section*{Acknowledgement}
This research of Y.K. is partially supported by NSF-DMS (1313312). The research of K.M is partially supported by  the  Department of Education GAANN (P200A120192).  We would like to give a special appreciation to Marisabel Rodriguez and Daniel Burkow for their assistances. \\}


{\small
\bibliography{PP_dispersal2.bib}
}


\end{document}